\documentclass[11pt]{article}

\usepackage[T1]{fontenc}
\usepackage[latin1]{inputenc}

\usepackage{amsfonts}
\usepackage{amssymb}
\usepackage{amsmath}
\usepackage[english,french]{babel}
\usepackage[all]{xy}
\usepackage{latexsym}
\usepackage{epsfig,epic,eepic}
\usepackage{psfrag}
\usepackage{eufrak}
\newcommand{\goth}[1]{\EuFrak{#1}}

\newcommand{\so}{\goth{so}}

\newcommand{\dd}{\,\text{\rm d}}

\newcommand{\FF}{{\cal F}}
\newcommand{\GG}{{\cal G}}

\newcommand{\R}{\mathbb R}

\newcommand{\T}{\mathbb T}
\newcommand{\Z}{\mathbb Z}

\newcommand{\Id}{\operatorname{Id}}

\newcommand{\Diff}{\operatorname{Diff}}

\newcommand{\Fr}{\operatorname{Lor}}
\newcommand{\saut}{\vspace{0.2cm}}
\newcommand{\SO}{\mbox{\rm SO}}
\newcommand{\oO}{\mbox{\rm O}}
\newcommand{\GL}{\mbox{\rm GL}}
\newcommand{\SL}{\mbox{\rm SL}}
\newcommand{\gl}{\mbox{\rm gl}}
\newcommand{\Sbb}{\mbox{$\mathbb S$}}

\newcommand{\RR}{\mathbb R}

\usepackage{theorem}
\newtheorem{enonce}{}[section]
\newtheorem{thm}[enonce]{Theorem}
\newtheorem{cor}[enonce]{Corollary}
\newtheorem{prop}[enonce]{Proposition}
\newtheorem{lem}[enonce]{Lemma}
\newtheorem{defi}[enonce]{Definition}
\newtheorem{fact}[enonce]{Fact}
\newtheorem{claim}[enonce]{Claim}
\theorembodyfont{\rmfamily}
\newtheorem{rema}[enonce]{Remark}

\newenvironment{proof}{\begin{trivlist}\item {\noindent{\bf Proof.}}}{\hfill $\Box$\end{trivlist}}

\psfrag{(A-,A+)=q}{$(A_-,A_+)=q$}
\psfrag{(A+,A-)=p}{$(A_+,A_-)=p$}
\psfrag{(A+,A+)}{$(A_+,A_+)$}
\psfrag{(A-,A-)}{$(A_-,A_-)$}
\psfrag{B1}{$B_1\times{\mathbb S}^1$}
\psfrag{B2}{$B_2\times{\mathbb S}^1$}
\psfrag{B3}{$B_3\times{\mathbb S}^1$}
\psfrag{P(B1)}{$P(B_1)$}
\psfrag{P(B2)}{$P(B_2)$}
\psfrag{P(B3)}{$P(B_3)$}
\psfrag{fig. 1}{{\bf Figure 1.}}
\psfrag{fig. 2}{\phantom{!}\begin{tabular}[t]{c}{\bf Figure 2.} The $P(B_i)$ are given for\\ correspondence with Figure 3\end{tabular}}
\psfrag{identification}{identification}
\psfrag{circle}{\phantom{!}\parbox[t]{6cm}{Circle $\{(\frac{\pi}{2},\frac{\pi}{2})\}\times{\mathbb S}^1$, along which\\both lightlike foliations have\\an horizontal tangent space.}}
\title{Lorentzian foliations on $3$-manifolds}
\author{C.\@ Boubel, P.\@ Mounoud, C.\@ Tarquini}
\date{\today}

\begin{document}
\selectlanguage{english}
\maketitle
\begin{abstract}
We study transversely Lorentzian foliations on the closed $3$-manifolds.
We classify them under a completeness hypothesis and we deduce the dual classification of codimension $1$ geodesically complete timelike totally geodesic foliations. Besides we provide an example of a Lorentzian foliation on a compact $3$-manifold which is equicontinuous on a {\em proper} dense open subset of the manifold, and hence not transversely complete.
\end{abstract}
\selectlanguage{french}
\begin{abstract}
Nous \'etudions les feuilletages transversalement lorentziens sur les vari\'et\'es ferm\'ees de dimension $3$. Sous une hypoth\`ese de compl\'etude, nous les classifions. Nous en d\'eduisons la classification, duale, des feuilletages de codimension $1$ de type temps totalement g\'eod\'esiques et g\'eod\'esiquement complets. Par ailleurs nous fournissons un exemple de feuilletage lorentzien sur une $3$-vari\'et\'e compacte, equicontinu sur un ouvert dense {\em propre} et donc non transversalement complet.
\end{abstract}
\selectlanguage{english}
\vskip .5cm 
\noindent
{\small \sl Keywords: Lorentzian foliation, transversely complete, timelike totally geodesic foliation}\\
{\sl \small Mathematical subject classification: 53C12; 53C50}

\noindent\hrulefill
\section{Introduction}
In this article we study the foliations admitting a transverse  Lorentzian metric. It is the  Lorentzian analogue of the better known Riemannian foliations and they still have a rigid transverse geometry. We recall that a transverse  (pseudo\nobreakdash-)Riemannian metric of a foliation $\FF$ is a (pseudo\nobreakdash-)Riemannian metric on the normal bundle $\nu(\FF)=TM/T\FF$, invariant under the flow of any vector field tangent to $\FF$.
We send the reader to the book \cite{Molino} of P.\@ Molino for a general theory in the Riemannian case and to section \ref{section_definitions} here for a precise definition. The Lorentzian foliations are thus the foliations admitting a transverse pseudo-Riemannian metric of index $1$.

One of the main differences between pseudo-Riemannian (and in particular Lorentzian) and Riemannian geometry is that a pseudo-Riemannian metric may admit non-equiconti\-nuous sequences of isometries. In the same way, the 
holonomy pseudogroup of a transversely pseudo-Riemannian foliation may be non equicontinuous.
For example, on the compact $3$-manifolds, \'E.\@ Ghys used that the Anosov flows having smooth weak stable and unstable foliations and preserving a volume form are transversely Lorentzian to classify them, see \cite{Ghys-fourier}.

Nevertheless, on compact manifolds, we may hope that equicontinuity fails in very special situations.  Indeed, so is it for sequences of isometries of a compact pseudo-Riemannian manifold: such a sequence, either is equicontinuous, or has nowhere bounded derivative, see \cite{Zeghib-jdg}. Moreover, an equicontinuous group of pseudo-Riemannian isometries preserves a Riemannian metric, this is mainly the Ascoli theorem. The Lorentzian foliations that were already known satisfied also such an alternative ``everywhere or nowhere equicontinuous'': they were either Anosov or Riemannian.

However, it turns out that this alternative is false for Lorentzian foliations. In section \ref{noncompletesection} we give an example of a Lorentzian foliation which is equicontinuous on a dense proper open subset. However, adding a "completeness" hypothesis, called transverse completeness, see Definition 
\ref{complete} here, this alternative is still true. Using Molino's theory and the classification of the Lorentzian Anosov flows, 
given by \'E.\@ Ghys in \cite{Ghys-fourier}, we prove the following theorem:\saut

\noindent{\bf Theorem \ref{alg-ano/riem}} {\slshape Up to finite cover, a $1$-dimensional transversely complete Lorentzian foliation on a compact $3$-manifold is either smoothly equivalent to the foliation generated by an algebraic Anosov flow or Riemannian.}\saut

\noindent Using then Y.\@ Carri\`ere's classification \cite{Car} of the Riemannian flows on compact $3$-manifolds, we give a detailed panorama of the last case.

As a natural corollary of Theorem \ref{alg-ano/riem}, we classify the codimension $1$ foliations of Lorentzian compact $3$-manifolds, the leaves of which are \emph{timelike}, totally geodesic and geodesically complete for a Lorentzian  metric. Indeed the integral foliation of the orthogonal distribution of those foliations are transversely complete Lorentzian foliations, see Propositions \ref{geodcomplete=>complete} and \ref{base}. We obtain:
\vskip .2cm
\noindent{\bf Theorem \ref{FTG}} {\slshape 
Up to a finite cover the timelike geodesically complete totally geodesic codimension $1$ foliations of the closed $3$ dimensional manifolds are
\begin{enumerate}
\item the foliations of the circle bundles over the torus, transverse to the fibres.
\item the foliations on $\T^3_A$ with a compact leaf and without Reeb nor type $II$ components.
\end{enumerate}
}
\vskip.2cm
\noindent
We send the reader to section  \ref{section classification 2-feuilletages} for the definition of a type $II$ component and to section \ref{completude} for a recall of that of $\T^3_A$. We just precise that, on $\T^3_A$, the foliations  without Reeb nor type $II$ components are exactly the taut foliations.

Finally to conclude this article we remark that, up to a finite cover, the tangent distributions of codimension $2$ Lorentzian foliations are given by the intersections of the kernels of two $1$-forms. We show then how to characterize those foliations in terms of $1$-forms.
\vskip .2cm \noindent
{\bf Structure of the article.} In section \ref{section_definitions} we introduce the notion of transversely complete pseudo-Riemannian foliation. In section \ref{section_alternative} we prove the alternative "Riemannian or Anosov" for transversely complete Lorentzian foliations and in section \ref{section_classification} we give the classification: Theorem \ref{alg-ano/riem} and its detailed version Theorem \ref{cas Riemannien}.
 
 Section \ref{section classification 2-feuilletages} deals with timelike totally geodesic foliations and Theorem \ref{FTG} is proved. 
 In  section \ref{noncompletesection}, we give examples of non transversely complete Lorentzian foliations and finally we give in section \ref{formeslineaires} the properties of the one-forms defining a Lorentzian codimension $2$ foliation.\\

\noindent{\bf Smothness of the objects.} All manifolds involved here are without boundary and all differential objects are of class $C^2$.\\

\noindent
{\bf Thanks.} The idea of this work was given us by a weekly workshop about conformally Lorentzian foliations at the \'Ecole Normale Sup\'erieure de Lyon in March-May 2003, together with Charles Frances, Thierry Barbot and Abdelghani Zeghib. We thank them for their contributions to this workshop, where we all learnt good mathematics.

\section{Definitions of completeness}\label{section_definitions}
\subsection{The definitions}
Let $\FF$ be a foliation of codimension $q$ on a manifold $M$. We will say that $\FF$ is transversely Lorentzian if it admits a so-called transverse Lorentzian metric. One can define this transverse structure using a holonomy pseudogroup $H$ of $\FF$ acting on some transverse manifold $T$ and a Lorentzian metric $g_T$ on $T$ such that each element of $H$ preserves $g_T$. However we thought that a definition in terms of a bundle-like metric would be clear enough. We recall what are a bundle-like metric, a basic function, form or  vector field, and what is a transverse parallelism, for a foliation $\FF$ on a manifold $M$ (see for instance \cite{Molino}).
\begin{enumerate}
\item  A real function of $M$ is basic (or foliated) for $\FF$ if it is constant along the leaves.
\item A vector field $Y$ on $M$ is basic (or foliated) for $\FF$ if for all vector field $X$ tangent to $\FF$ the vector field $[X,Y]$ is tangent to $\FF$. This is equivalent to ask that the flow of $Y$ preserves the foliation.
\item A $1$-form $\alpha$ on $M$ is basic (or foliated) for $\FF$ if for all vector field $X$ tangent to $\FF$ we have $\alpha(X)=0$ and $i_X\dd\alpha=0$ (or equivalently if $\alpha(X)=0$ and $L_X\alpha=0$ for $L_X$ the Lie derivative of $X$). Notice that this constraint is stronger than that of point 2 for vectorfields. Indeed, a foliation $\FF$ always admits non-zero basic vectorfields ---yet little interesting---: the fields tangent to $\FF$, while a foliation may admit no non-zero basic form.
\item A transverse parallelism of $\FF$ is a $q$-tuple of basic vector fields $(X_1,\dots,X_q)$ spanning a complement of $T\FF$ at every point of $M$. In general, a foliation $\FF$ does not admit any parallelism. If it admits one, it is said to be transversely parallelizable.
\item A (pseudo-)Riemannian metric $g$ on $M$ will be called $\FF$-bundle-like if the restriction of $g$ to $T\FF$, the tangent bundle of $\FF$, is non-degenerate and if
$$
 L_X \left(g_{|T\FF^\perp}\right)=0,
$$ for any vector field $X$ tangent to $\FF$.
\end{enumerate}

Notice that the integral curves of a nowhere lightlike and nowhere vanishing Killing vector field gives a particular case of foliation having a bundle-like metric.

\begin{defi}
We will say that $\FF$ is transversely Lorentzian if there exists a $\FF$-bundle-like  metric $g$  which induces on the normal bundle $\nu(\FF)=TM/T\FF$ a Lorentzian metric.
\end{defi}

We denote by $\GL(M,\FF)$ the transverse frame bundle of $\FF$, it is a $\GL(q)$-principal bundle. We can lift the foliation $\FF$ on a foliation $\FF^1$ on $\GL(M,\FF)$ which has the same dimension but is of course of codimension $q+q^2$. Saying that $\FF$ is transversely Lorentzian is equivalent to saying that $\GL(M,\FF)$ has a $\FF^1$-foliated reduction to the Lorentzian group $\oO(1,q-1)$ i.e. a reduction to the group $\oO(1,q-1)$ of $\GL(M,\FF)$ which is $\FF^1$-saturated. We will denote by $\Fr(\FF)$ this reduction.

\begin{claim}\label{claim}
The transverse Lorentzian structure of $\FF$ gives a natural transverse parallelism of $\FF^1$ on $\GL(M,\FF)$, hence on $\Fr(\FF)$.
\end{claim}

This is a consequence of the existence of the Levi Civita connection of the Lorentzian metric and of the application of the theory of Molino for Riemannian foliation \cite{Molino}. Following Molino we introduce the:

\begin{defi}
Let $\GG$ be a transversely parallelizable foliation of the manifold $N$ and let $L_c(\GG)$ be the set of complete vector fields of $N$ basic for the foliation $\GG$. We will say that $\GG$ is a {\em complete} transversely parallelizable foliation if for all points $x$ in $N$ the set $L_c(\GG)(x)=\{X(x), X\in L_c(\GG)\}$ spans a complement to $T\GG$.
\end{defi}

\begin{defi}\label{complete}
We will say that a Lorentzian foliation $\FF$ is transversely complete if the foliation $\FF^1$ on $\Fr(\FF)$ is a complete transversely parallelizable foliation.
\end{defi}

\begin{rema}\label{GL ou Fr}
{\rm
The last definition can be equivalently stated with $\GL(M,\FF)$ instead of $\Fr(\FF)$. Indeed for all points $z\in\GL(M,\FF)$ and all $h\in\GL(q)$ let us denote by $R_h(z)=z.h$ the right action of $\GL(q)$ on $\GL(M,\FF)$. Then for all vector $Z$ in $\gl(q)$, the Lie algebra of $\GL(q)$, we have a natural complete vectorfield on $\GL(M,\FF)$ tangent to the fibres of $\GL(M,\FF)\to M$ given by
$$
\widetilde Z(z)=\left.\frac{\dd}{\dd t}\right|_{t=0} z.\exp(tZ)
$$
This gives, for all elements $z$ in $\GL(M,\FF)$, a canonical identification of $\gl(q)$ to the tangent space of the fibre of $\GL(M,\FF)\to M$ at $z$. Hence the transverse completeness of a Lorentzian foliation depends on the existence of complete vectorfields in $L_c(\FF^1)$, transverse to $\FF^1$ and to those fibres.}
\end{rema}

This remark implies also that Riemannian foliations are transversely complete since in this case we have a $\FF^1$-foliated reduction of $\GL(M,\FF)$ to the group $\oO(q)$ which is compact.
A well known fact for transversely parallelizable complete foliations is the following proposition (see \cite{Molino}):

\begin{prop}\label{Diff-transitif}
Let $(N,\GG)$ be a transversely parallelizable complete foliation on a connected manifold $N$. Then the group $\Diff(N,\GG)$ of diffeomorphisms of $N$ keeping $\GG$ invariant (i.e.\@ all elements of $\Diff(N)$ sendind leaves of $\GG$ onto leaves of $\GG$), acts transitively on $N$.
\end{prop}

The proof follows from the fact that the complete parallelism ensures that the orbits of $\Diff(N,\GG)$ are open subsets of $N$. Thus there is only one orbit.

It is well known for Riemannian $\FF$-bundle-like metrics that a geodesic starting orthogonally to $T\FF$ remains orthogonal to it for any time. This is still true in the pseudo-Riemannian setting. Thus it is meaningful to speak about the geodesic completeness of the orthogonal distribution $T\FF^\perp$.

\begin{defi}\label{geodesicallycomplete}
A Lorentzian foliation $\FF$ will be called geodesically complete if there exists a bundle-like Lorentzian metric for which the orthogonal distribution $T\FF^\perp$ is geodesically complete.
\end{defi}

It gives another notion of completeness. Be careful that the existence an $\FF$-bundle-like Lorentzian geodesically complete metric does not imply that all $\FF$-bundle-like Lorentzian metrics are geodesically complete (see Proposition \ref{pticorollaire} in Section \ref{section classification 2-feuilletages}). The following proposition connects those definitions of completeness.

\begin{prop}\label{geodcomplete=>complete}
A geodesically complete Lorentzian foliation is transversely complete in the sense of Definition \ref{complete}.
\end{prop}

\begin{proof}
We suppose that $\FF$ is a geodesically complete Lorentzian foliation of codimension $q$ and let $\pi:\GL(M,\FF)\to M$ be the principal bundle of transverse frames to $\FF$.

We first construct the geometric connection on $\GL(M,\FF)$ corresponding to the Levi-Civita covariant derivative on $M$ given by the bundle-like Lorentzian metric. For all $z\in\GL(M,\FF)$, $z$ is a transverse frame at $x=\pi(z)$ that is to say $z$ is a linear isomorphism from $\RR^q$ to $\nu_x(\FF)$. Let $X$ be a tangent vector to $M$ at $x$ orthogonal to $\FF$. We can find a $C^1$ curve $\gamma:]-\varepsilon,\varepsilon[\to M$ orthogonal to $\FF$ and such that $\gamma(0)=x$ and $\gamma'(0)=X$. We denote by $\tau_\gamma$ the parallel transport along $\gamma$. Thus for each $t$ in $]-\varepsilon,\varepsilon[$, the linear isomorphism $\tau_\gamma(t)\circ z$ is a transverse frame over $\gamma(t)$. Now we can set out 
$$
s_{x,z}(X)=\left.\frac{\dd}{\dd t}\right|_{t=0} \tau_\gamma(t)\circ z
$$
It is well defined since $\tau_\gamma$ depends smoothly on $t$ and its $1$-jet depends only on $x, X$ and $z$. We set $H_z$ to be $s_{x,z}(T\FF^\perp)$. From the definition it satisfies $H_{z.h}=\dd_zR_h(H_z)$ for all $h\in\GL_q(\RR)$ (where $R_h$ is the right translation by $h$). And we also have $\dd_z\pi(s_{x,z}(X))=X$. Therefore $H_z$ is right invariant and transverse to the fibres and to the foliation $\FF^1$ (the lifting of $\FF$ on $\GL(M,\FF)$). We can define a geometric connection $P_z$ at $z$ by setting $P_z=H_z\oplus T_z\FF^1$.

Let $\theta$ be the fundamental form on $\GL(M,\FF)$. We recall its definition: for all $z\in\GL(M,\FF)$ and all vector $Z$ tangent at $z$ to $\GL(M,\FF)$, we have $\theta_z(Z)=z^{-1}(\overline{\dd_z\pi(Z)})$ where $X\to\overline X$ denotes the projection from $TM$ to $\nu(\FF)$. The $1$-form $\theta$ is $\RR^q$-valued, basic for $\FF^1$ and moreover satisfies $R^*_h\theta=h^{-1}\theta$ (see \cite{Molino}). Restricted to $H_z$, the fundamental form $\theta_z$ defines a linear isomorphism onto $\RR^q$. Let $(e_1,\dots,e_q)$ be the canonical basis of $\RR^q$, for all $i$ we denote by $Z_i(z)$ the unique vector of $H_z$ such that $\theta_z(Z_i(z))=e_i$. Then
\begin{itemize}
\item for all $i$, the vector field $Z_i$ if basic for $\FF^1$ since so is $\theta$.
\item $(Z_1(z),\dots,Z_q(z))$ is a basis of $H_z$ and projected to $TM$ by $\dd_z\pi$, it gives a basis of $T\FF_x^\perp$.
\end{itemize}

By Remark \ref{GL ou Fr} to finish the proof we just have to show that each $Z_i$ is complete. So let $z_i(t)$ be the integral curve of $Z_i$ starting at $z_i(0)=z$. It is defined on an open neighborhood $I$ of $0$ in $\RR$. We denote by $\gamma_i(t)=\pi(z_i(t))$ its projection on $M$. Then $\gamma_i$ is orthogonal to $\FF$. We will prove that it is a geodesic. Indeed recall that:
$$
\nabla_{\gamma_i'}\gamma_i'(x)=\lim_{t\to 0}\frac 1 t \left(\tau_{\gamma_i}(t)^{-1}(\gamma_i'(t))-\gamma_i'(0)\right)
$$
By definition for all $t$ in $I$ we have
$$
z^{-1}(\overline{\gamma_i'(0)})=\theta_z(Z_i(0))=e_i=\theta_{z_i(t)}(Z_i(z_i(t)))=\theta_{z_i(t)}(z_i'(t))=z_i(t)^{-1}(\overline{\gamma_i'(t)})  
$$

\noindent Moreover if $\tau_i(t)$ denotes the curve $\tau_{\gamma_i}(t)\circ z$, since $\pi(\tau_i(t))=\gamma_i(t)=\pi(z_i(t))$ we have $\dd_{z_i(t)}\pi(\tau_i'(t))=\dd_{z_i(t)}\pi(z_i'(t))$ and as $\tau_i'(t)$ and $z_i'(t)$ belong to $H_{z_i(t)}$ we have the equalities $\tau_i'(t)=z_i'(t)$ and $z_i(t)=\tau_{\gamma_i}(t)\circ z$. Thus $\tau_{\gamma_i}(t)^{-1}(\gamma_i'(t))=\gamma_i'(0)$ and $\nabla_{\gamma_i'}\gamma_i'(x)=0$. This argument shows that $\nabla_{\gamma_i'}\gamma_i'$ vanishes on $I$ therefore $\gamma_i$ is a geodesic orthogonal to $\FF$. As by hypothesis $\gamma_i$ is complete the curve $z_i$ is too and then $Z_i$ is a complete basic vector field.
\end{proof}

The reciprocal of this proposition is not true: there exist non geodesically complete but transversely complete Lorentzian foliations. 

\subsection{An example of a transversely complete but not geodesically complete Lorentzian foliation}\label{completude}

\begin{prop}
Let $\FF$ be the strongly stable (or unstable) foliation of the Anosov flow obtained by suspending a hyperbolic diffeomorphism of the 2-torus. Then $\FF$ is transversely Lorentzian complete in the sense of Definition \ref{complete}, and not geodesically complete in the sense of Definition \ref{geodesicallycomplete}.
\end{prop}

\begin{proof} The manifold obtained by this suspension is classically denoted by ${\mathbb T}^3_A$, with $A\in{\rm SL}_2(\Z)$ the linear hyperbolic automorphism of ${\mathbb T}^2$ to which the initial diffeomorphism is isotopic. It is well known that $\FF$ is a transversely affine Lie foliation, see \cite{Car}. Hence it is transversely parallelizable and therefore both Riemannian and Lorentzian. It follows from remark \ref{GL ou Fr} that $\FF$ is transversely complete.  We are going to prove that $\FF$ admits no geodesically complete bundle-like Lorentzian metric.

Let $\widetilde \FF$ be the lift of $\FF$ to the universal cover of ${\mathbb T}^3_A$. The space of leaves of $\widetilde \FF$ is diffeomorphic to $\R^2$. The fundamental group of ${\mathbb T}^3_A$, denoted by $\Gamma$, acts on this space. This action can be seen as generated by the transformations:
$$(x,y)\mapsto (x+\tau,y)\quad\text{and:}\quad
(x,y)\mapsto (\lambda x, y+1),$$
where $\tau$ belongs to a \emph{dense} subgroup of $\R$ and $\lambda$ is an irrational number greater than $1$.
A transverse metric for $\FF$ is just a metric on $\R^2$ invariant by this action. We are going to prove that such metrics are never geodesically complete. As a first example we can notice that this action can be seen as the action of a subgroup of the group of isometries of a half Minkowski plane (with lightlike boundary) on it. This is certainly not a geodesically complete metric.

Unfortunately the space of transverse metrics of $\FF$ is big, the reason is that it is transversely parallelizable and it admits basic functions. (See an example of a non-homogegeous transverse Lorentzian metric in section \ref{formeslineaires} p.\@ \pageref{metriquenonplate}.)
The strategy of the proof is to find lightlike incomplete geodesics on the space on leaves using the works of Y.\@ Carri\`ere and L.\@ Rozoy. They proved in \cite{Ca-Ro} that the completeness of lightlike geodesics on Lorentzian surfaces is related to the holonomy of the lightlike foliations. More precisely we will use the fact that almost all lightlike geodesics accumulating on a closed lightlike geodesic are incomplete.

Let $g$ be a Lorentzian metric on $\R^2$ invariant by the action of $\Gamma$. This metric has two transverse lightlike geodesic foliations $\mathcal G^1$ and $\mathcal G^2$. Those foliations are invariant under $\Gamma$.

 Let us first suppose that one of the lightlike foliations of $g$, for example $\mathcal G^1$, is everywhere transverse to the foliation $\mathcal G^0$  the leaves of which are the $\R\times\{y\}$. Note that this is in particular the case if one of those foliations is $\mathcal G^0$. This foliation is determined by one of its leaves, the other being deduced by the translations  $(t,0)$, for any $t\in \R$ . It means that there exists a real number $\alpha$ such that for any $x\in \R$ the leaf of $\mathcal G^1$ through $(x,0)$ contains the  point $(x+\alpha,1)$.
Consequently the leaf containing $(\frac{\alpha}{\lambda - 1},0)$ is fixed by the transformation $\gamma : (x,y)\mapsto (\lambda x, y+1)$.
Let us consider $\langle\gamma\rangle$ the group generated by $\gamma$. It acts isometrically freely and properly on $\R^2$. Thus $\R^2/ \langle\gamma\rangle$ is a cylinder endowed with the quotient metric $\bar g$.
The foliation $\mathcal G^1$ gives a lightlike foliation $\overline{\mathcal G^1}$ which has a closed attractive leaf. By \cite{Ca-Ro}, $\bar g$, and therefore $g$, is not geodesically complete.

We suppose now that $\mathcal G^1$ is somewhere but not everywhere tangent to $\mathcal G^0$. It implies that $\mathcal G^1$ has at least a countable number of leaves in common with $\mathcal G^0$. Let $\gamma'$ be a translation belonging to $\Gamma$. Repeating the same argument as above with $\gamma'$ instead of $\gamma$, we see that $g$ must be geodesically incomplete.
\vskip .2cm
We have proved hat $\R^2$ has no $\Gamma$-invariant geodesically complete Lorentzian metric and therefore that $\FF$ is not a geodesically complete Lorentzian foliation.
\end{proof}

\section{The alternative: Riemannian or Anosov}\label{section_alternative}

\begin{thm}\label{alternative}
Let $\FF$ be a $1$ dimensional complete transversely Lorentzian foliation on a compact connected manifold $M$ of dimension three. Then either $\FF$ is Riemannian (and thus transversely parallelizable) or up to a finite cover (in fact a $2$-, $4$- or $8$-cover) the foliation $\FF$ is given by the orbits of an Anosov flow.
\end{thm}

\begin{proof}
Up to a $2$-cover, $\FF$ can be supposed to be orientable, and then parametrized by some flow $\varphi^t$. Both isotropic directions $E_1$ and $E_2$ of the Lorentzian metric span each fibre of the normal bundle $\nu(\FF)=TM/T\FF$. They are necessarily preserved or swapped by the differential of the flow. Up to an additional possible $2$ or $4$-cover (as $\oO(1,1)$ has four connected components) we also assume that the bundles $E_1$ and $E_2$ are {\em each one} preserved by the flow, and preserved with an orientation i.e. they admit non vanishing sections $s_1$ and $s_2$. We define the {\it deformation cocycle} $u(x,t)$ by:
$$
\varphi^t_*(s_1(x))=e^{u(x,t)}s_1(\varphi^t(x))\qquad\mbox{and:}\qquad
\varphi^t_*(s_2(x))=e^{-u(x,t)}s_2(\varphi^t(x)),
$$
where $\varphi_*^t$ denotes the action of $\varphi^t$ on $\nu(\FF)$, and where $u$ satisfies:
$$u(x,s+t)=u(x,t)+u(\varphi^t(x),s).$$

The leaves of $\FF^1$ are given by the action of $\varphi^t$ on $\Fr(\FF)$. More precisely the action of $\varphi^t$ on $M$ gives rise to an action on $\nu(\FF)$ and then to an action on $\GL(M,\FF)$ which stabilizes $\Fr(\FF)$. We will denote by $\varphi_*^t$ this action.

The foliation $\FF^1$ on $\Fr(\FF)$ is by definition a complete parallelizable foliation. Therefore by Proposition \ref{Diff-transitif}, either all leaf of $\FF^1$ is precompact in $\Fr(\FF)$ or no leaf is.
\medskip

\noindent{\bf In the first case} we apply Theorem 2 of R. Wolak in \cite{Wolak} (in this case the precompacity of leaves of $\FF^1$ implies the completeness see \cite{ptitnote}), it shows that the foliation $\FF$ is Riemannian. Normalizing $s_1$ and $s_2$ by a bundle-like Riemannian metric we have a transverse parallelism for $\FF$. 
\medskip

\noindent{\bf In the second case} we will show that the flow $\varphi^t$ is Anosov with $E_1$ and $E_2$ (the projections on $\nu(\FF)$ of) its weak stable and unstable bundles. The idea is to show that there exists a closed leaf of $\FF^1$ in $\Fr(\FF)$. Using Proposition \ref{Diff-transitif} this will imply that all leaf of $\FF^1$ is closed and using the structure Theorem of Molino (see Theorem \ref{theoreme molino} below) the flow $\varphi$ will be quasi-Anosov hence Anosov by a result of Ma\~ne.
\medskip

We will say that a flow $\varphi^t$ on a compact $C^\infty$ manifold $M$ is quasi-Anosov over a saturated subset $K$ (or that $K$ is quasi-Anosov for $\varphi^t$) if for all vector $v\in \nu(\FF)$ based on a point of $K$, $v\neq 0$, the set $\{\|\varphi^t_*(v)\|, t\in\RR\}$ is unbounded. Here $\varphi_*$ is the flow on $\nu(\FF)$ induced by the differential of $\varphi^t$ and the foliation $\FF$ is given by the orbits of $\varphi$. We will use the following proposition which is a weak version of a Theorem of Ma\~ne (\cite{Mane}).

\begin{prop}\label{proposition Mane}
Let $K$ be a saturated connected compact subset of $M$ such that $\varphi$ is quasi-Anosov over $K$. Then $K$ is hyperbolic for $\varphi$.
\end{prop}

\noindent As T.\@ Barbot said us, the arguments of \cite{Mane} are remarkably simplified here. So they are worth to be reproduced to give the following, elementary proof. We thank T.\@ Barbot for it.\saut

\begin{proof}
By assumption, no $\varphi^t_\ast$-orbit is bounded. Therefore, for any $x$ in $K$, there exists a real number $T(x)$ such that $u(x,T(x))$ is bigger than $1$ (apply $\varphi^t_\ast$ to $s_1(x)$). Applying similary $\varphi^t_\ast$ to $s_2(x)$, we find a $T'(x)$ such that $u(x,T'(x))$ is smaller than $-1$. By continuity of $u$ and compactness of $K$ we deduce the existence of some positive real number $T$ such that:
$$
\forall x\in K,\ \sup_{|t|\leqslant T} u(x,t) \geqslant 1\quad\mbox{and}\quad
\inf_{|t|\leqslant T} u(x,t) \leqslant -1
$$

Let us set $C=\sup_{0\leqslant t\leqslant T, x\in K} |u(x,t)|$.

\begin{lem}\label{encadrement}
For any real positive numbers $0\leqslant s\leqslant t$, and any point $x$ of $K$, we have:
$$
\min(-C,-C+u(x,t))\leqslant u(x,s)\leqslant\max(C,u(x,t)+C)
$$
\end{lem}

\begin{proof} Fix $t>0$, $x\in K$ and denote by $\tau$ a time between $0$ and $t$ such that $u(x,\tau)$ is equal to $\sup_{0\leqslant s\leqslant t} u(x,s)$. By definition of $T$, the segment $[\tau-T,\tau+T]$ cannot be contained in $[0,t]$. Hence either $0\leqslant\tau< T$ in which case we have $u(x,s)\leqslant u(x,\tau)\leqslant C$ for every $0\leqslant s\leqslant t$, or $0\leqslant t-\tau< T$, in which case $u(x,s)\leqslant u(x,\tau)=u(x,t)-u(\varphi^{\tau}(x),t-\tau)\leqslant u(x,t)+C$.

\noindent The proof of the lower bound is similar.
\end{proof}

\noindent At every point $x$ of $K$, we obtain the following dichotomy:\\
- either $u(x,t)$ tends to $+\infty$ when $t$ tends to $+\infty$,\\
- or there is some increasing sequence $t_n$ converging to $+\infty$ for which $u(x,t_n)$ is uniformly bounded from above. By Lemma \ref{encadrement}, in this case all the $u(x,s)$, for $0\leqslant s$, are bounded from above.

Similarly, the $u(x,s)$, for $0\leqslant s$, are bounded from below or converge to $-\infty$. Symmetrically, we obtain similar statements for $t$ converging to $-\infty$.

Since the norm of $e^{u(x,t)}s_1(x)$ cannot be bounded from above for all times $t$, $u(x,t)$ must tend to $+\infty$ for $t$ converging to $+\infty$ or $-\infty$. Furthermore, in the first case, $u(x,t)$ tends to $-\infty$ for negative $t$, and in the last case, $u(x,t)$ tends to $+\infty$ for negative $t$ (Hint: consider $e^{-u(x,t)}s_{2}(x)$).

In other words, there is an $\varepsilon(x)=\pm 1$ such that $\lim_{t\to+\infty}\varepsilon(x)u(x,t)=+\infty$. By compactness of $K$, there is a positive real $T'$ such that $\varepsilon(x)u(x,T')\geqslant 2$ for every $x$. It follows, by connexity of $K$, that $\varepsilon(x)$ does not depend on $x$, let us say that $\varepsilon(x)=1$ everywhere.

\begin{lem}\label{final Mane}
If there exists two reals $T'$ and $a$ such that $\inf_{x\in K} u(x,T')\geqslant a>0$, then $K$ is hyperbolic for $\varphi^t$.
\end{lem}

\begin{proof}
Possibly reversing the time we can suppose that $T'$ is positive.

\noindent Let $C=\sup_{0\leqslant t\leqslant T', x\in K} |u(x,t)|$, it is a finite constant. For any positive $t$ and any point $x\in K$, considering the unique integer $n$ such that $nT'\leqslant t <(n+1)T'$, we observe that:
\begin{align*}
u(x,t)&=u(x,T')+u(\varphi^{T'}(x),T')+\cdots+u(\varphi^{nT'}(x),t-nT')\\
&\geqslant a(n-1)-C\\
&\geqslant a\frac t {T'}-C\\
&\geqslant a't+b' \quad\mbox{for some constants $a'$, $b'$ with $a'>0$.}
\end{align*}

Then at $x$ vectors in $E_1$ are exponentially expanded in the future and exponentially contracted in the past. The fact that $\varphi$ is transversely Lorentzian proves that, at $x$, elements of $E_2$ are exponentially contracted in the future and are exponentially expanded in the past.\medskip

It is well known that such a behavior implies the Anosov property, but since all readers are not supposed to be familiar with this notion, we provide here the line of the proof.

Let $E_{11}$, respectively $E_{22}$, the pull back of $E_1$, $E_2$, under the projection $TM\to\nu(\FF)$. The exponential behavior of $\varphi_{\ast}^{t}$ on $E_1$, $E_2$ implies the hyperbolic property as soon as there are, over $K$, some $D\varphi^{t}$-invariant $1$-dimensional subbundle $E^{uu}$, respectively $E^{ss}$, of $E_{11}$, respectively $E_{22}$ transverse to ${\mathcal F}$. To exhibit such invariant subbundles, for example $E^{uu}$, identify the space ${\mathcal E}$ of $1$-dimensional subbundles of $E_{11}$ complementary of $\langle X\rangle=T\FF$ ($X$ is the vector field of the flow), with the space of continuous functions $f$ on $K$ {\it via\/} some splitting $E_{11} = \langle X\rangle \oplus \langle \eta\rangle$, where $\eta$ is some vector field: $f$ defines a subbundle by $E_{f}(x) = \langle f(x)X + \eta \rangle$.

The action of $D\varphi^{t}$ on $E_{11}$ is defined by functions $u(x,t)$ and $\alpha(x,t)$ such that:
\begin{align*}
 D\varphi^{t}(x, \eta(x)) = & e^{u(x,t)}\eta(\varphi^{t}(x)) + \alpha(x,t)X(\varphi^{t}(x)) \\
 D\varphi^{t}(x, X(x)) = & X(\varphi^{t}(x))
\end{align*}

The natural action induced by $D\varphi^{t}$ on $\mathcal E$ defines a $1$-parameter group of transformations ${\mathcal A}_{t}$ of $C^0(K)$ given by:
\[ {\mathcal A}_{t}(f)(x) = 
e^{-u(\varphi^{-t}(x),t)}(f(\varphi^{-t}(x))+\alpha(\varphi^{-t}(x), t)) \]
Since $u(x,T) \geq a>0$, the operator ${\mathcal A}_{T}$ is contracting on $C^0(K)$ equipped with the sup-norm. It admits a unique fixed point $f_{0}$. Since all the ${\mathcal A}_{t}$ commute with ${\mathcal A}_{T}$, $f_{0}$ must be a common fixed point for all the ${\mathcal A}_{t}$. Then, $E_{f_{0}}$ is then the required invariant bundle.
\end{proof}
This lemma concludes the proof of Proposition \ref{proposition Mane}.
\end{proof}

To prove the following corollary and for the end of the proof we need the:

\begin{thm}[P. Molino \cite{Molino} Theorem 4.2]\label{theoreme molino}
Let $\GG$ be a complete transversely parallelizable foliation on a connected manifold $N$. Then $N$ is the total space of a fibre bundle $p:N\to W$ on some manifold $W$ and  the closure of the leaves of $\GG$ are the fibres of $p$. Moreover, on each fibre of $p$, $\GG$ is a transversely Lie foliation with dense leaves.
\end{thm}

\begin{cor}\label{corollaire Mane}
If $\varphi$ is quasi-Anosov over a saturated connected compact subset $K$ of $M$ or if there exists a closed non compact leaf in $\Fr(\FF)$, then $\varphi$ is Anosov.
\end{cor}

\begin{proof}
We show that if $\varphi^t$ is quasi-Anosov on such a $K$, then $\Fr(\FF)$ has a closed non compact leaf. Then we prove that this implies that $\varphi^t$ is Anosov.

By Proposition \ref{proposition Mane} the compact set $K$ is hyperbolic. Therefore every leaf $F^1$ in $\pi^{-1}(K)$ is closed in $\Fr(\FF)$ and non compact (recall that $\pi:\Fr(\FF)\to M$). Indeed $F^1$, parametrized by $\varphi^t_\ast$, tends to infinity when $t$ tends to $\pm\infty$. This means by Theorem \ref{theoreme molino} applied to $(\Fr(\FF),\FF^1)$, that all leaves of $\FF^1$ are the fibres of a fibre bundle $p:\Fr(\FF)\to W$. As the foliation is given by a flow, this fibre bundle becomes a $\R$-principal bundle.

The application $s:M\to\Fr(\FF),x\mapsto(s_1(x),s_2(x))$ is a non basical section of the principal bundle $\pi:\Fr(\FF)\to M$. Using the continuity of the action of $\RR$ (by $\varphi^t_\ast$) on $\Fr(\FF)$, the triviality of the fibre bundle $p:\Fr(\FF)\to W$ (since $\RR$ is contractible) and the compacity of $s(M)$, we prove that there exists a real number $T>0$ such that $\varphi^T_*(s(M))\cap s(M)=\emptyset$. Then by definition of $T$, $\inf_{x\in M} |u(x,T)|=a>0$. By continuity we can suppose that $a=\inf_{x\in M} u(x,T)$. Applying Lemma \ref{final Mane} with $K=M$ we see that $\varphi^t$ is Anosov.
\end{proof}

To finish we have to prove the existence of either a saturated quasi-Anosov connected compact subset $K$ of $M$ or a closed non compact leaf of $\FF^1$, the corollary \ref{corollaire Mane} will conclude.

If there exists a periodic orbit of $\varphi^t$ in $M$, say $L$, then there exists a time $T>0$ such that for all $x\in L$ we have $\varphi^T(x)=x$ and by the cocycle relation $u(x,nT)=nu(x,T)$ for all integer $n$. Since by assumption $\{u(x,t), t\in\R\}$ is unbounded, $u(x,T)\neq 0$ and $\varphi^t$ is quasi-Anosov over $L$. So we can suppose that no orbit of $\varphi^t$ is periodic.

Let $L^1$ be a leaf of $\FF^1$ and suppose that the closure of $L^1$ is tangent to a fibre of $\pi$ at a point $z$, hence there exists a vector $Z$ tangent to $\overline{L^1}$ at $z$ and such that $\dd_z\pi(Z)=0$. We can extend this vector by a vector field $\widetilde{Z}$ tangent to the fibre and basic for $\FF^1$ (see Remark \ref{GL ou Fr}). By minimality of $\overline{L^1}$ for $\FF^1$, the vector field $\widetilde{Z}$ is tangent to $\overline{L^1}$ at every point. Integrating this vector field by a right translation along the fibres of $\pi$ we see that $\overline{L^1}$ is saturated by the fibres of $\pi$, i.e. $\overline{L^1}=\pi^{-1}(\pi(\overline{L^1}))$ or equivalently $\pi(\Fr(\FF)\setminus \overline{L^1})=M\setminus \pi(\overline{L^1})$. Since the projection is open, the set $\pi(\overline{L^1})$ is closed hence equal to $\overline{\pi(L^1)}$. Then $\pi(\overline{L^1})$ is a connected compact $\FF$-saturated subset of $M$. By Theorem \ref{theoreme molino}, any leaf of $\overline{L^1}$ is dense in $\overline{L^1}$ hence for all $x$ in $\pi(\overline{L^1})$ the set $\{u(x,t), t\in\R\}$ cannot be bounded from above or from below. This means that $\pi(\overline{L^1})$ is quasi-Anosov for $\varphi^t$.

Otherwise the closure of $L^1$ is everywhere transverse to the fibres of $\pi$, we can suppose that it is true for all leaf of $\FF^1$. There are three cases according to the dimension of $\overline{L^1}$.

If $\dim\overline{L^1}=1$ the leaf $L^1$ is closed and non compact, we are done.

If $\dim\overline{L^1}=2$, then we use the fact that $\FF^1$ restricted to if $\overline{L^1}$ is a $G$-Lie foliation where $G$ is a simply connected Lie group of dimension $1$ hence equal to $\R$. The developing map is complete, i.e. it is a locally trivial fibre bundle, this comes form the completeness of the transverse parallelism of $\FF^1$. As moreover the leaves are simply connected, the holonomy morphism $\rho$ from $\pi_1(\overline{L^1})$ to $\R$ is injective. This implies that $\pi_1(\overline{L^1})$ is abelian. Now since $\overline{L^1}$ is a non compact surface, $\pi_1(\overline{L^1})$ is a free group (this comes from the fact that a non compact surface has a proper Morse function, bounded from below and without local maximum see \cite{phillips}). Therefore $\pi_1(\overline{L^1})$ is equal to $\Z$. But as $\overline{L^1}$ is minimal, $\rho(\pi_1(\overline{L^1}))\cong\Z$ must be dense in $\R$. This is impossible.

If $\dim\overline{L^1}=3$, as $\overline{L^1}$ is transverse to the fibres, $\pi(\overline{L^1})$ is open in $M$. But this holds for all closures of leaves of $\FF^1$. Hence $\pi(\overline{L^1})$ is open and closed in $M$, which is connected, so $M=\pi(\overline{L^1})$. As $\overline{L^1}$ is non compact, $\pi|_{\overline{L^1}}$ is not injective, i.e. there is a point $z\in\overline{L^1}$ and a $g\in\SO_0(1,1)$, $g\neq \Id$, such that $z.g$ belongs to $\overline{L^1}$. This implies that, for all integer $n$, the equality $\overline{L^1}g^n=\overline{L^1}$. So, as before, by minimality for all $x$ in $M=\pi(\overline{L^1})$ the set $\{u(x,t), t\in\R\}$ is bounded neither from above nor from below. This precisely saying that the flow $\varphi^t$ is quasi-Anosov and Anosov (by Proposition \ref{proposition Mane}). But this implies that all leaves of $\FF^1$ are closed and non compact. Finally this case is impossible too.
\end{proof}

\section{Classification}\label{section_classification}

\subsection{Anosov flows preserving a volume form are algebraic, after Ghys}\label{section Ghys}

In \cite{Ghys-fourier}, \'E. Ghys proves that, on compact $3$-manifolds and up to conjugation by a $C^\infty$-diffeomorphism, the Anosov flows the weak stable and unstable distributions of which are $C^\infty$ and which preserve a volume form (notice that those properties does not depend on their parametrization, so are properties of the foliation) are exactly the {\em ``algebraic Anosov flows''} i.e.\@ the finite coverings and the finite quotients of the following:\medskip

$\bullet$ the geodesic flows of the unit tangent bundle of the hyperbolic compact surfaces,\medskip

$\bullet$ the suspensions of the linear hyperbolic diffeomorphisms $A$ of the $2$-torus (those give rise to a manifold denoted by ${\mathbb T}^3_A$).\medskip

\noindent Notice that both are transversely Lorentzian ---the first ones with constant nonvanishing transverse curvature (positive or negative, it amounts to the same as the signature of the metric is $(1,1)$, ``neutral''), the second ones with a flat Lorentzian metric.

\begin{thm}\label{alg-ano/riem}
Up to finite cover, a $1$-dimensional transversely complete Lorentzian foliation on a compact $3$-manifold is either smoothly equivalent to an algebraic Anosov flow or Riemannian.
\end{thm}

\begin{proof}
The only thing we still have to prove is the completeness of an algebraic Anosov flow in the sense of definition \ref{complete}. This follows by the proposition \ref{geodcomplete=>complete}.
\end{proof}

\subsection{With Carri\`ere's work: the detailed panorama}\label{section Carriere}

Finally, using Carri\`ere's classification, provided in \cite{Car}, of the Riemannian flows on the compact $3$-manifolds, we detail the description of the foliations of Theorem \ref{alternative}, in the case they are at once Lorentzian and Riemannian. This provides this last statement:

\begin{thm}\label{cas Riemannien} Up to conjugation by a diffeomorphism, the $1$-dimensional transversely  complete Loren\-tzian foliations of the compact $3$-manifolds are:

\noindent{\bf (a)} Up to a possible $2$-, $4$- or $8$-covering, the algebraic Anosov flows on ${\mathbb T}^3_A$ and on the unit tangent bundle of the compact hyperbolic surfaces,

\noindent{\bf (b)} The simultaneously Riemannian and Lorentzian ones. Up to a possible $2$- or $4$-covering, those are the parallelizable foliations of the compact $3$-manifolds, which are, up to a possible orientation $2$-covering:
\begin{trivlist}
\item[{\bf \hspace*{0.5cm}(i)}] the linear flows on a 3-torus,
\item[{\bf \hspace*{0.5cm}(ii)}] the flows defined by the fibres of a circle bundle over the 2-torus,
\item[{\bf \hspace*{0.5cm}(iii)}] the flows defined by the stable (resp. unstable) direction of the suspension of a linear hyperbolic 
diffeomorphism of the 2-torus.\saut
\end{trivlist}
\end{thm}

\noindent Of course, categories {\bf (b)\,(i)} and {\bf (b)\,(ii)} have a nonempty intersection.
\vskip .2cm

\begin{proof}
We first recall that a transversely Lorentzian and Riemannian foliation on a compact manifold is transversely complete Lorentzian (see remark \ref{GL ou Fr}).
Moreover, up to a possible $2$- or $4$-covering, foliations of codimension $2$  and at once Riemannian and Lorentzian are parallelizable.
Indeed, the bundle of transverse frames which are Riemannian orthonormal and Lorentzian orthogonal is a covering with fibre $\SO(2,\R)\cap\oO(1,1)\simeq (\Z/2\Z)^2$. As in claim \ref{claim}, the lift of $\FF$ to this bundle is transversely parallelizable.
Conversely, a parallelizable foliation preserves many pseudo-Riemannian or Riemannian metrics: take a transverse parallelism and say it is a (pseudo-)orthonormal basis.

Now we are left with classifying the parallelizable $1$-dimensional foliations. Taking again, if necessary, a 2-cover we can assume the foliation to be orientable. A dimension $1$ orientable foliation is usually called a flow yet without considering any parametrization. Instead of simply using Y.\@ Carri\`ere's classification of the Riemannian flows on the compact $3$-manifolds (see \cite{Car} part III A),  we pick up the relevant parts of his proof to achieve the desired list. Here is Carri\`ere's ground theorem, from which follows his classification.

\begin{thm}[Y.\@ Carri\`ere, \cite{Car} II.C, Theorems 3+2]\label{theoreme carriere}
The closures of the leaves of a Riemannian flow $\FF$  are tori, on which
${\cal F}$ induces a dense linear flow.
\end{thm}

\noindent{\bf Remark.} In the very special case of a parallelizable flow, in which we are interested, this result follows immediately from (a more complete version of) the Molino theorem \ref{theoreme molino} and from a result of Carri\`ere's thesis (see \cite{Car}, part II.C, Theorem 1).

So let $\FF$ be a parallelizable flow on a compact $3$-manifold. After the Molino Theorem \ref{theoreme molino}, the closure of its leaves have a common dimension and are fibres of a fibration on some compact manifold $W$. After Theorem \ref{theoreme carriere}, these fibres are tori and the study splits into three cases according to their dimension.\saut

1. There is a dense leaf of ${\cal F}$. Then $\FF$ is a dense linear flow on the 3-torus. Conversely, such a flow is immediately parallelizable. This is the case {\bf (b) (i)} of Theorem \ref{cas Riemannien}.\saut

2. {All the leaves of $\FF$ are closed. Then $M$ is a circle bundle over a closed surface $W$}. Since  ${\cal F}$ is transversely Lorentzian, $W$ admits some Lorentzian metric and thus it must be a 2-torus ($\FF$ is transversely orientable). In this case, the leaves of ${\cal F}$ are defined by the fibres of a locally trivial fibration by circles over the 2-torus. Conversely, the foliation defined by such a fibration is automatically parallelizable. This is the case {\bf (b) (ii)} of Theorem \ref{cas Riemannien}.\saut

3. {The closure of any leaf of $\FF$ is a $2$-torus. Then $M$ is a torus bundle over a compact $1$-manifold, thus over the circle.} Such a bundle is obtained by considering the product $\T^2 \times [0,1]$ and gluing the tori $\T^2 \times \{ 0 \}$ and $\T^2 \times \{ 1 \}$ by some diffeomorphism $\varphi$ of $\T^2$. On each fibre $\T^2 \times \{ s \}$, the flow ${\cal F}$ induces a linear flow of irrational direction $\Delta_s$. But since $\Delta_s$ must be continuous in $s$, this direction is actually constant;  we denote it by $\Delta$. Thus $\varphi$ lets the direction $\Delta$ invariant. By Lemma I.B.5 of \cite{Ghys-Serj}, $\varphi$ is in fact isotopic to a linear diffeomorphism $A$ of $\T^2$, through an isotopy $\varphi_t$ such that for all $t \in [0,1]$, $\varphi_t$ preserves the direction $\Delta$. Our flow $\FF$ is then smoothly conjugated to the flow induced by the translation of direction $\Delta$ on the manifold $\T_A^3$. Now a matrix $A$ of $\SL_2(\Z)$ preserving an irrational direction is either the identity or an hyperbolic matrix.

If $A=\Id$ we are again in the case {\bf (b) (i)} of theorem \ref{cas Riemannien}. Else, $\Delta$ is the stable or the unstable direction of $A$ and we are in the case {\bf (b) (iii)}. It remains to check that any flow of that last type is parallelizable. In fact, it is even transversely a Lie flow, modeled on the affine group $\mbox{\rm AG}(1,\R)$ (see for instance \cite{Car} I.D. example 6). Thus any left-invariant frame on $\mbox{\rm AG}(1,\R)$ yields a transversely invariant framefield.
\end{proof}
\subsection{Isometric foliations}
Among transversely Lorentzian $1$-dimensional foliations take place, in particular, those given by  nowhere vanishing spacelike Killing fields: the invariant metric of the manifold gives an invariant transverse metric on $(T\FF)^\perp$, canonically identified with $\nu(\FF)$. So let us define:

\begin{defi}
A transversely (pseudo\nobreakdash-)Riemannian flow is called isometric if it admits a Killing parametrization for some bundle-like metric on the manifold.
\end{defi}

We have the following natural characterization, see  \cite{Car} Proposition III.B.1.

\begin{prop}
A transversely pseudo-Riemannian flow is isometric if and only if it admits a parametrization which preserves a codimension $1$ transverse distribution.
\end{prop}

Moreover an algebraic Anosov flow has smooth strongly stable and unstable distributions. Hence it preserves the distribution spanned by them and it is isometric.
This proposition shows also that a flow which is both Lorentzian and Riemannian is Lorentzian isometric if and only if it is Riemannian isometric. Thus, thanks to \cite {Car} corollary III.B.4, we see that  the only non isometric transversely complete Lorentzian flows on compact $3$-manifolds are  of type {\bf (b)\,(iii)}. Actually the nowhere vanishing spacelike Killing fields of $3$-dimensional manifolds have been  classified by A. Zeghib in \cite{Zeghib-jdg}. It turns out that they  all generate transversely complete Lorentzian foliations.
It means that the examples given in section \ref{noncompletesection} are not isometric.

\section{The geodesically complete totally geodesic timelike foliations in dimension 3}\label{section classification 2-feuilletages}
Let us first specify what is a totally geodesible foliation of a pseudo-Riemannian manifold.
\begin{defi}
Let $\mathcal G$ be a codimension $1$ foliation on a manifold $M$.
The foliation $\mathcal G$ is said to be  geodesible if there exists a pseudo-Riemannian metric for which 
the leaves of $\mathcal G$ are totally geodesic submanifolds, i.e.\@ any geodesic of the metric starting tangentially to a leaf stays in the leaf.
\end{defi}
In Riemannian geometry it is well known that a transversely oriented codimension $1$ foliation is geodesible if and only if it is 
transverse to a Riemannian flow. It is the starting point of the classification given by Y.\@ Carri\`ere and \'E.\@ Ghys, in 
\cite{Carriere-Ghys}  for $3$-dimensional manifolds and in \cite{Ghys-totgeod} in any dimension.
We will see that the  Lorentzian analog is less general. 

The signature of the restriction of the metric to some totally geodesic, connected submanifold $N$ is constant, as any tangent space $T_xN$ is stable by parallel transport along any path drawn on $N$. Thus, there exist three types of totally geodesic connected hypersurfaces in a Lorentzian manifold. They can be timelike, spacelike, or lightlike if this restriction is respectively Lorentzian, Riemannian, or degenerate. Hence we can talk about spacelike, timelike and lightlike totally geodesic foliations if all their leaves are of this type. (In general however, there may be leaves of each type, see \cite{ken} and \cite{mounoud2}). 
The spacelike case is exactly the same as when $M$ is Riemannian and is classified. The lightlike one has been principally studied by A.\@ Zeghib in \cite{Zeghib} (see also \cite{mounoud2}) but a classification is only conjectured. In the general case  a classification is out of reach for the moment (cf. \cite{mounoud2}). Here we will look at the timelike case and to begin this study we specify the definition.
\begin{defi}
Let $\mathcal G$ be a codimension $1$ foliation on a manifold $M$.
The foliation $\mathcal G$ is said to be timelike geodesible if there exists a Lorentzian metric $g$ for which 
the leaves of $\mathcal G$ are totally geodesic, timelike submanifolds.
\end{defi}
We state the following fact, the proof of which is exactly the same as the Riemannian one (see \cite{Carriere-Ghys}).
\begin {prop}\label{base}
A smooth transversely oriented codimension $1$ foliation is timelike (resp.\@ spacelike) totally geodesible if and only if it is transverse to a transversely Lorentzian (resp.\@ Riemannian) flow.
\end{prop}
Or more precisely:
\begin {prop}\label{base2}
Let $(M,g)$ be a Lorentzian manifold. 
A smooth codimension $1$ timelike (or spacelike) foliation $\mathcal G$ is totally geodesic if and only if $g$ is bundlelike for the foliation generated by its orthogonal distribution.
\end{prop}
As we said, we restrict ourselves to smooth \emph{timelike} totally geodesic foliations. They are not the only ones but they still represent an important and interesting family. Moreover they are definitely related to Lorentzian flow. We further restrict to the geodesically complete case i.e.\@ when the geodesics tangent to the foliation are all complete. Thanks to proposition \ref{geodcomplete=>complete} and \ref{base} and theorem \ref{cas Riemannien} we can state:
\begin{cor}\label{triv}
Let $\GG$ be a totally geodesic complete codimension $1$ foliation on a closed $3$-manifold $M$. Then, up to finite cover, one of the following situations occurs: 
\begin{enumerate}
\item The manifold $M$ is a circle bundle over the torus and $\mathcal G$ is transverse to the fibres.
\item The foliation $\mathcal G$ is transverse to the strong (un)stable direction of a flow defined by the suspension of a linear hyperbolic diffeomorphism of the $2$-torus.
\item The foliation $\mathcal G$ is transverse to an algebraic Anosov flow.
\end{enumerate}
\end{cor}

\begin{proof}
By Proposition \ref{base} our foliations are transverse to a Lorentzian flow, moreover this flow is necessarily transversely complete by proposition \ref{geodcomplete=>complete}. So Theorem \ref{cas Riemannien} gives more or less the list.
The only thing is to notice that if $\mathcal G$ is transverse to a translation flow on a torus with irrational slope it is also transverse to a translation flow with rational slope.
\end{proof}

There exist foliations satisfying 1.\@ 2.\@ or 3.\@; they are actually all timelike geodesible but it is not true that they can all be realized as \emph{geodesically complete} timelike totally geodesic foliations.
The first and second cases correspond to both spacelike and timelike totally geodesible foliations.  Thus they have been studied by Y.\@ Carri\`ere and \'E.\@ Ghys in \cite {Carriere-Ghys}. In the first case the foliations are given by the suspensions of pairs of commuting diffeomorphisms of the circle. We will see that they can be realized as geodesically complete timelike totally geodesic foliations. We will not try to classify them.
 In the second case, Y. Carri\`ere and \'E.\@ Ghys showed that the foliations have no compact leaves and therefore, using \cite{Ghys-Serj}, are diffeomorphically conjugated to the weak stable or unstable foliation of the Anosov flow. 
We have seen in subsection \ref{completude} that even if they are  transversely complete, the strongly stable and unstable foliations of the Anosov flow obtained by suspension do not have any geodesically complete bundle-like metrics. This construction will not give geodesically complete timelike totally geodesic foliations.
\vskip.1cm
We are going to describe the foliations transverse to an algebraic Anosov flow i.e.\@ either to the geodesic flow of a hyperbolic surface $\Sigma$, and the ambient manifold is the unitary tangent bundle $T_1\Sigma$, or the flow given by the suspension of an hyperbolic diffeomorphism  of the torus $A$, and the ambient manifold is the torus bundle over the circle with monodromy $A$ classically denoted by $\T^3_A$. Of course those flows are just considered as oriented $1$ dimensional foliations.
In what follows $\FF$ will denote one of the algebraic Anosov flows and $\mathcal G$ will denote a codimension $1$ foliation transverse to $\FF$.
\\
We are going to prove the following proposition :
\begin{prop}\label{trans}
A codimension $1$ foliation $\mathcal G$ is transverse to an algebraic Anosov flow if and only if it is one of the following.
\begin{enumerate}
\item It is conjugated to the weak stable or unstable foliation of an algebraic Anosov flow.
\item $M=\T^3_A$ and $\mathcal G$ has a compact leaf and is without Reeb nor type $II$ components.
\end{enumerate}
\end{prop}

To be able to give the definition of a type II component, we have  to recall the construction of the two fundamental classes of examples of codimension $1$ foliations on $\T^2\times [0,1]$ tangent to the boundary, (cf.\@ the work of R.\@ Moussu and R.\@ Roussarie \cite{rrww}). 

We start with a foliation on the annulus $\Sbb^1\times [0,1]$, invariant by rotation, tangent to the boundary and without compact leaves in the interior. There are two possibilities: this foliation is given by the suspension of a diffeomorphism of the interval without fixed point in the interior, or it is a Reeb component. Now let us consider the product foliation on the annulus times $[0,1]$ and at last let us glue both sides of this product after some rotation. We obtain this way two families of foliations on $\mathbb \T^2\times [0,1]$ tangent to the boundary and without compact leaves in the interior. According to the choice of the rotation the noncompact leaves are all cylinders or all planes.
\\
Those foliations can also be obtained by the following 1-forms on $\T^2\times [0,1]$
$$\Omega_1=dr+\varphi(r)\omega\quad \mathrm{and}\quad \Omega_2=(1-2r)dr+\varphi(r)\omega$$
where $r\in [0,1]$, $\omega$ is a  linear form on the torus and $\varphi$ is a smooth function such that $\varphi(0)=\varphi(1)=0$ and $\varphi(r)>0$ if $0<r<1$. Consequently we denote them $\mathcal G([\omega],1)$ and $\mathcal G([\omega],2)$.

We define another class of foliations. We start from the foliation $\mathcal G([\omega_0],2)$ where the form $\omega_0$ defines a trivial foliation by circles on $\T^2$  (i.e.\@ where the rotation is the identity and of course the leaves of $\mathcal G([\omega_0],2)$ are cylinders). There exists an orientation preserving diffeomorphism of order $2$, without fixed points which permutes the connected components of the boundary and preserves the foliation. When we take the quotient we obtain a foliation on a manifold with boundary $\T^2$ foliated by cylinders accumulating on the boundary. In fact the manifold is the non trivial fibre bundle over the Klein bottle by $[0,1]$ (see \cite{rouss} for details). Such a foliation will be called a cylindrical component.

Now we can set the following definition.
\begin{defi} 
A foliation will be called a component of type $II$ if it is topologically conjugated to a foliation $\mathcal G([\omega],2)$ (i.e.\@ obtained from a Reeb component on an annulus) or to a cylindrical component.
\end{defi}

This definition is very close from the definition p. 104 of \cite{rouss} but we do not ask the foliations to be with cylinder leaves. 

\begin{defi}
An embedded surface $S$ of a $3$-manifold $M$, not diffeomorphic to a $2$-sphere, will be called incompressible if the map $\pi_1(S)\rightarrow \pi_1(M)$ is injective.
\end{defi}

We recall that the Novikov theorem says that a codimension $1$ foliation on a compact $3$-manifold has a compressible leaf if and only if it contains a Reeb component. Now we can start the proof of proposition \ref{trans}.
\medskip

\noindent{\bf Proof of proposition \ref{trans}.} Let us consider a foliation $\mathcal G$ transverse to an algebraic Anosov flow and prove that it is in case 1.\@ or 2.\@ We have seen in \ref{section Ghys} that a transverse metric preserved by an algebraic Anosov flow $\FF$ has constant curvature --- it comes from the fact that $\FF$ has a dense orbit. If $\FF$ is a suspension flow the transverse metric has to be flat and if $\FF$ is a geodesic flow the metric is with non zero constant curvature (as seen in \ref{section Ghys}, its sign has no meaning on a $2$ dimensional Lorentzian manifold).

Moreover if $L$ is a compact leaf of $\mathcal G$, it is Lorentzian, so has to be a torus or a Klein bottle.
But by the Gauss-Bonnet-Avez theorem (see \cite{Avez-Gauss-Bonnet}) it cannot be endowed with a metric with non vanishing curvature. Hence a foliation transverse to the geodesic flow has no compact leaves. Consequently we can use the work of \'E.\@ Ghys \cite{Ghys-anosov} and see that $\mathcal G$ is then conjugated to the weak stable foliation of $\FF$. 

We still have to consider the suspension case. If $\mathcal G$ has no compact leaves then according to \cite{Ghys-Serj} it is conjugated to either the weak stable, or unstable, foliation of $\FF$.
 Suppose now that the foliation has a compact leaf $L$.  
\begin{fact}\label{nono}
The foliation $\mathcal G$ contains neither Reeb nor type $II$ components.
\end{fact}

\begin{proof}
The first point is to know if a compact leaf $L$ can separate the manifold, i.e.\@ is $\T^3_A \setminus L$ connected ?

\begin{lem}\label{fa}
A compact leaf $L$ of $\mathcal G$ cannot separate $\T^3_A$.
\end{lem}

\begin{proof} 
Suppose that $L$ separates the manifold into two manifolds $M_1$ and $M_2$. As the Anosov flow is a flow, thus oriented, all its orbits cutting $L$ go, for example, from $M_1$ to $M_2$. Thus any orbit of $\FF$ cutting $L$ is trapped in $M_2$. This is in contradiction with the fact that an algebraic Anosov flow possesses a dense set of closed orbits.
\end{proof}

\noindent Lemma \ref{fa} implies that $\mathcal G$ has no Reeb nor cylindrical components. We still have to prove that $\mathcal G$  does not contain some component topologically conjugated to some $\mathcal G([\omega],2)$. Taking again a possible finite cover we assume the Anosov flow transversely oriented.
As the  orientation of the compact leaves of such a foliation are opposite, such a component admits no transverse path joining its two compact leaves. As an algebraic Anosov flow has a dense leaf it cannot be transverse to it.
\end{proof}

\begin{rema}\label{taut}{\rm
Fact \ref{nono} actually just says that the foliation $\mathcal G$ is \emph{taut}, i.e.\@ every compact leaf admits a closed transversal.
Moreover a  foliations is taut if and only if there exists a metric (a priori Riemannian) for which the leaves are minimal surfaces (see \cite{sullivan}). Hence Riemannian totally geodesic foliations are taut but it can be verify that timelike totally geodesic foliations are  as well taut and fact \ref{nono} is true for these foliations.}\end{rema}

We are done with the direct part of Proposition \ref{trans}. Let us now consider a foliation $\mathcal G$ as in point 1.\@ or 2.\@ of  Proposition \ref{trans}. Taking a possible $2$ cover we assume that both $\FF$ and $\mathcal G$ are oriented.

If $\mathcal G$ is as in point 1., we notice that the pull-back of the Anosov flow  by the flow of a strong unstable (resp.\@ stable) vector field is transverse to it and of course still Anosov.

If $\mathcal G$ is as in point 2., it is a foliation on $\T^3_A$ without Reeb nor cylindrical components but with a compact leaf; it follows from theorem 2.1 of \cite{good} that the compact leaf can not separate the manifold. Hence if we cut along it we obtain a foliation on a  manifold with incompressible ---by the Novikov theorem and as $\mathcal G$ contains no Reeb component---, non connected boundary , i.e.\@  $\T^2\times [0,1]$, according to theorem 4-2 of \cite{EM}.

Using the fact that the leaves are incompressible tori and reproducing the beginning of the proof of Theorem 2 in \cite{rrww}, we see that $\mathcal G$ is isotopic to a foliation $\mathcal G'$ whose compact leaves are  fibres of the fibration of $\T^3_A$ on $\Sbb^1$. We use now Theorem 1 of \cite{rrww} in the following form.
\begin{thm}[R.\@ Moussu and R.\@ Roussarie]\label{roberts}
Let $\mathcal G$ be a foliation on $\T^2\times [0,1]$ tangent to the boundary and without compact leaves in the interior. 
Then there exists a homeomorphism of $\T^2\times [0,1]$ fixing the boundary 
which sends the foliation $\mathcal G$ 
to a model foliation $\mathcal G([\omega],1)$ or $\mathcal G([\omega],2)$.
\end{thm}
If follows from this theorem that there exists a homeomorphism of $\T^3_A$ which sends $\mathcal G$ transversely to the natural Anosov flow of $\T^3_A$. As being transverse is an open property and as for $3$-dimensional manifolds (see \cite{munkres}) the set of diffeomorphisms is dense inside the set of homeo\-morphisms, we see that there is actually a diffeomorphism sending $\mathcal G$ transversely to the Anosov flow. This completes the proof of Proposition \ref{trans}.\hfill$\Box$
\vskip .2cm\noindent
We are going to prove now that the stable or unstable foliation $\GG$ of an algebraic Anosov flow is never timelike totally geodesic and geodesically complete. According to Corollary \ref{triv} and the result of subsection \ref{completude}, the only case left we have to deal with is when $\mathcal G^\perp$ is an algebraic Anosov foliation we denote by $\FF$. As the weak (un)stable foliations of two algebraic Anosov foliations on some manifold $M$ are conjugated, we take $\varphi$ a diffeomorphism of $M$ such that $\GG=\varphi(\GG')$ with $\GG'$ the weak (un)stable foliation of $\FF$. The foliation $\FF$ has a unique transverse Lorentzian $(G,X)$ structure. Indeed, the lightlike foliations of $\FF$ are its weak stable and unstable foliations: this determines the conformal class of the transverse Lorentzian metric. Moreover $\FF$ has dense leaves, so this metric is unique up to a scalar factor and has constant curvature. The transverse model space $X$ is then the Minkowski plane or the universal cover of the two dimensional de Sitter space. It means that the lift of $\FF$ to the universal cover $\widetilde M$ of $M$ is given by the preimages of a submersion $\delta$ from $\widetilde M$ to $X$ called the developing map. As $\FF$ is conjugated to an algebraic Anosov flow and as its transverse $(G,X)$ structure is unique, this structure is conjugated to that of the algebraic flow, which is complete (see \cite{Ghys-fourier}), that is to say, $\delta$ is a locally trivial fibration with fibre $\R$.

The restriction of $\delta$ to the lift $\widetilde L$ of a leaf of $\mathcal G$ gives it a $(G,X)$-structure.
If $\mathcal G$ is geodesically complete then the restriction of $\delta$ of $\FF$ to $\widetilde L$ is a diffeomorphism. It means that the couple $(\FF, \mathcal G)$ must lift to the universal cover as a product foliation. We are going to prove that this last point is not possible.

As $\GG=\varphi(\GG')$ with $\GG'$ the weak (un)stable foliation of $\FF$, $\GG$ has a cylindrical leaf $L_0=\varphi(L'_0)$. The lift $\widetilde{L'_0}$ of $L'_0$ to $\widetilde M$ is a plane containing the lift of a compact leaf of $\FF$, so is invariant by an element $\gamma\in\pi_1(M)$ sent by the holonomy representation $\rho$ to an isometry $\rho(\gamma)$ of $X$ having a fixed point. Now the lift $\widetilde{L_0}$ of $L_0$ on $\widetilde M$ is fixed by $\widetilde{\varphi}^{-1}\circ\gamma\circ\widetilde{\varphi}=\gamma$, as the lift $\widetilde{\varphi}$ of $\varphi$ to $\widetilde{M}$ commutes with the action of $\pi_1(M)$. Eventually, as $\delta_{|\widetilde{L_0}}:\widetilde{L_0}\rightarrow X$ conjugates the action of $\gamma$ and of $\rho(\gamma)$, $\gamma\in\pi_1(M)$ fixes a point in $\widetilde{L_0}\subset\widetilde{M}$, which is impossible. It means that $\mathcal G$ can not be timelike totally geodesic and geodesically complete.
\vskip.2cm
The last point is to prove that the other cases are geodesically complete. This is done by noticing that, in those cases, the couple $(\mathcal G, \FF)$ lifts to the universal cover as a product and that  $\FF$ possesses a Lorentzian and geodesically complete transverse $(G,X)$-structure. 

\begin{rema} Actually those two conditions are necessary and sufficient. The main step to prove it is the following theorem proved in \cite {Po-Ri} by R.\@ Ponge and and H.\@ Rieckzigel.
\begin{thm}[R.\@ Ponge and H.\@ Rieckzigel]
If a timelike totally geodesic codimension $1$ foliation $\mathcal G$ is geodesically complete then the pair of foliations $\mathcal G$ and $\mathcal G^\perp$ lifts to the universal cover as a product.
\end{thm} 
In particular this theorem says that the orthogonal foliation of a geodesically complete timelike codimension $1$ foliation  possesses a \emph{complete} Lorentzian $(G,X)$-transverse structure. But to have a reciprocal we have to ask the Lorentzian manifold $X$ to be geodesically complete.
\end{rema}

Setting together what we just showed with Corollary \ref{triv} and Proposition \ref{trans} we obtain the following theorem:
\vskip .2cm
\begin{thm} \label{FTG}
Up to a finite cover the geodesically complete timelike  geodesible codimension $1$ foliations of the closed $3$ dimensional manifolds are:
\begin{enumerate}
\item The foliations of the circle bundles over the torus, transverse to the fibres.
\item The foliations on $\T^3_A$ with a compact leaf but without Reeb nor type $II$ components.
\end{enumerate}
\end{thm}

By the way we showed also the following little result, announced after Definition \ref{geodesicallycomplete}.

\begin{prop}\label{pticorollaire}
The algebraic Anosov flow $\FF$ of $\T^3_A$ admits some bundle-like Lorentzian metrics such that $T\FF^\perp$ is geodesically complete and some others such that it is not.
\end{prop}

\noindent
{\bf Remarks:}\medskip

1. In the second case the proof gives a description of the foliations. They contain a countable number of $\mathcal G([\omega],1)$  components (with possibly different $\omega$'s) and a closed set of toral leaves. Moreover, up to a finite cover, circle bundles over the torus are also torus bundles over the circle. Hence in the first case of theorem \ref{FTG}, the foliations having a compact leaf do have the same description but we have to consider also the foliations by planes and by cylinders.\medskip

  2.
 As we saw, a transversely Lorentzian flow $\FF$ on a $3$-manifold $M$ is always (after a possible $2$-cover) given by the intersection of two transverse codimension $1$ foliations $\mathcal G_1$ and $\mathcal G_2$.
Hence a Lorentzian totally geodesic foliation $\mathcal G$ together with those two foliations form what is called a total foliation: $(\mathcal G,\mathcal G_1,\mathcal G_2)$, i.e.\@ at any $p\in M$, $(T_p\mathcal G)^\perp\oplus (T_p\mathcal G_1)^\perp\oplus(T_p\mathcal G_2)^\perp=T^\ast_p M$. What about the converse ? Actually, we can give a Lorentzian geometric interpretation for $\mathcal G$ to be part of a total foliation of $M$. We recall that a diffeomorphism of a Lorentzian manifold $(M,g)$ is conform if and only if it preserves the lightcone of $g$. Therefore, up to a $2$-cover, we see that a dimension $1$ foliation on a $3$-manifold is transversely conformally Lorentzian if and only if it is given by the intersection of two transverse codimension $1$ foliations. As a codimension $1$ foliation is  umbilical (i.e.\@ its leaves are  umbilic for a metric of the manifold) if and only if it is transverse to a transversely conformal flow (adapt the argument given in \cite{Carriere-Ghys}, Prop I.3), we have that  a codimension $1$ foliation $\mathcal G$ is umbilical with timelike leaves for some Lorentzian metric if and only if there exist two codimension $1$ foliations $\mathcal G_1$ and $\mathcal G_2$ such that $(\mathcal G,\mathcal G_1,\mathcal G_2)$ is a total foliation of $M$. \medskip

 3.
There exists other timelike totally geodesic foliations on closed $3$-manifolds: The foliations transverse to Lorentzian flows which are not transversely complete. We will give such an example page \pageref{exampletotogeo}. It is not clear for the moment if it will be possible to classify them.  If we remove the hypothesis of completeness, the only thing we know is that the foliations are taut which is a weaker property. But it has interesting consequences. For example, as a taut foliation has no Reeb  components we know that the manifold is covered by $\R^3$. We can also notice that the taut foliations on  torus bundles over the circle are the foliations of corollary \ref{triv}. Hence they are all timelike totally geodesic. It means that we know all the timelike totally geodesic foliations of those manifolds.

\section{Non transversely complete Lorentzian flows}\label{noncompletesection}
\subsection{The examples}
We are going to give the simplest example we know of a non transversely complete Lorentzian flow. We will shortly say how to produce new examples from this one.
We first construct it as a transversely conformally Lorentzian flat foliation and show that it is Lorentzian. For a clearer exposition, we  replaced, when it was possible, technical details by drawings. A more precise and systematic studies will be done in a forecoming paper.

The torus $\R{\rm P}^1\times \R {\rm P}^1$ endowed with the metric $\dd \theta \dd \varphi$ is called the Einstein torus and will be denoted by Ein$_2$. Its main particularity is the size of its conformal group: it is a finite index extension of Diff$({\mathbb S}^1)\times$Diff$({\mathbb S}^1)$. This implies that being only a transversely conformally Lorentzian flat codimension $2$ foliation is too general for us.

But the topological group Diff$({\mathbb S}^1)\times$Diff$({\mathbb S}^1)$ contains PSL$(2,\R)\times$PSL$(2,\R)$ as a Lie subgroup. Thus we will rather construct our example among the foliations admitting a transverse $($PSL$(2,\R)\times$ PSL$(2,\R), \R {\rm P}^1\times \R {\rm P}^1)$-structure.
We will call them \emph{transversely Einstein} foliations. We recall that such foliations are defined by an equivariant submersion $D$, called the developing map, from  the universal cover of the ambient manifold $M$ to the Einstein torus. Moreover, this submersion must be equivariant ie there exists a morphism $\rho$ from the fundamental group of $M$ to PSL $(2,\R)\times$ PSL $(2,\R)$ such that $D\circ \gamma = \rho(\gamma)\circ D$ for any $\gamma$ in the fundamental group. Sometimes the universal cover is replaced by an intermediate covering called the holonomy cover, this will be the case here.
\vskip .2cm
Let us begin the construction. We denote by $\Sigma_2$ the 
compact orientable surface of genus two. In order to find a submersion from $\Sigma_2 \times \R$ to the Einstein torus, we first consider a projection $P=(P_1,P_2)$ from $\Sigma_2$ to Ein$_2$. We can easily find $P$ such that the image of $\Sigma_2$ is given by figure 1 and that any point in the interior of the image is regular (and has two preimages).
\begin{figure}[h]
\hspace*{\fill}
\epsfig{file=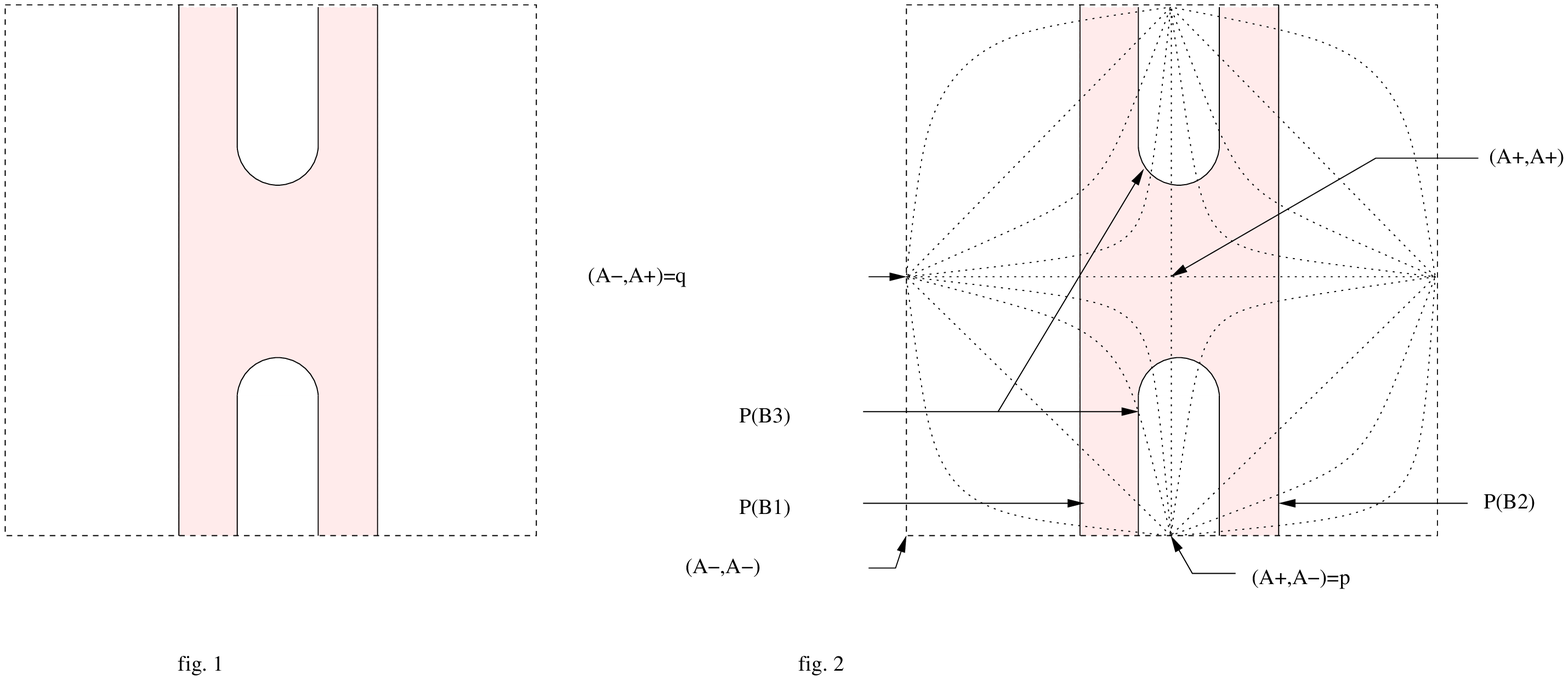,width=14cm}
\hspace*{\fill}
\end{figure}
We consider now a hyperbolic element $A$ of PSL$(2,\R)$. There exists an element $a$  of $\mathfrak{sl}(2,\R)$ such that $A= \exp a$.
In a well chosen parametrization of $\R$P$^1\times\R$P$^1$, and denoting by $A_+$ and $A_-$, respectively, the attractive and repulsive points af $A$, the orbits of $\exp (t a)\times \exp(-t a)$ on Ein$_2$ are the dotted curves of fig. 2.
We can now define our submersion: 
$$
\begin{array}{lccc}
& \Sigma_2 \times \R & \longrightarrow & \mathrm{Ein}_2\\
D:    &  (x,t) & \mapsto  & \left(\exp(t a)(P_1(x)), \exp(-t a)(P_2(x))\right)
\end{array}
$$
It is a submersion at any point $(x,t)$ with $P(x)$ in the interior of $P(\Sigma_2)$, the grey zone in Figure 1. It is also a submersion if $P(x)$ is in the one-dimensional manifold $\mathcal C$ forming the boundary of $P(\Sigma_2)$. Indeed, on the one hand, $\dd D(x,t).(T_x\Sigma_2\times\{0\})=T_{P(x)}{\mathcal C}$, and on the other hand, as the orbits of $\exp (t a)\times \exp(-t a)$ are transverse to $\mathcal C$, $\dd D(x,t).(\{0\}\times\R)$ is transverse to $T_{P(x)}{\mathcal C}$, this gives result.
The transformation $(x,t) \mapsto (x,t+1)$ is conjugated by $D$ to the transformation $(A,A^{-1})$ which lies in PSL$(2,\R)\times$PSL$(2,\R)$. Hence the fibres of $D$ go down to $\Sigma_2\times \R/\Z\simeq\Sigma_2\times{\mathbb S}^1$ and define on it a transversely Einstein foliation. 

The image of $D$, denoted by $U$, is Ein$_2$ deprived of the ``vertical'' circle $\{A_-\}\times\R{\rm P}^1$, containing the repulsive point $q=(A_-,A_+)$, and of the attractive point $p=(A_+,A_-)$. In order to prove that $D$ defines actually a Lorentzian foliation we have to find a metric on $U$ preserved by $(A,A^{-1})$. This will be done by finding  $U_A$ the biggest subset of Ein$_2$ on which $(A,A^{-1})$ acts isometrically and proving that $U$ is included in $U_A$. The points $p$ and $q$ are the attractive and repulsive fixed points of $(A,A^{-1})$ so they do not belong to $U_A$.
\begin{fact}
Let $A$ be a hyperbolic element of ${\rm PSL}(2,\R)$ and let $p$ and $q$ be the attractive and repulsive fixed points of $(A,A^{-1})$ acting on $\R {\rm P}^1\times\R {\rm P}^1$. There exists a Lorentzian metric (in the conformal class of the Einstein metric) on $\R {\rm P}^1\times\R {\rm P}^1 \setminus\{p,q\}$ preserved by $(A,A^{-1})$.
\end{fact}
\begin{proof}
Let us first consider the following map:
$$
\begin{array}{lccc}
& \R^2&\longrightarrow &\R/\pi\Z\times\R/\pi\Z\\
   \psi : &       (x,y) &\mapsto& (\arctan x, \arctan y)
\end{array}
$$
This map is a diffeomorphism on its image  and it provides a conformal embeding of the Minkowski $2$-space into the Einstein torus minus two lightlike circles. More precisely 
$$(\psi^{-1})^*\dd x\dd y = \frac{\dd \theta \dd \varphi}{\cos^2 \theta \cos^2 \varphi}.$$
Without any loss of generality we can suppose that there exists $\lambda\not\in\{0,1\}$ such that 
$$f_A(x,y) = \left( \psi^{-1}\circ (A,A^{-1})\circ \psi\right) (x,y) = (\lambda\, x,\frac{1}{\lambda}\, y)$$ 
and so leaves the metric $\dd x\dd y$ invariant. It is easy to see that $f_A$ leaves invariant any metric of the kind 
$\alpha (xy)\dd x\dd y$, where $\alpha :\R \rightarrow \R$ does not vanish.
In particular $f_A$ leaves invariant the metric $g_0=\frac{1}{1+x^2y^2} \dd x\dd y$, but 
$$(\psi^{-1})^*\frac{1}{1+x^2y^2} \dd x\dd y = \frac{\dd \theta \dd \varphi}{\cos^2 \theta \cos^2 \varphi+ \sin^2 \theta \sin^2 \varphi}.$$
This metric, a priori only defined on $\psi(\R^2)$, extends to $\R {\rm P}^1\times\R {\rm P}^1 \setminus\{p,q\}$. Moreover, as $\psi(\R^2)$ is dense in $\R {\rm P}^1\times\R {\rm P}^1 \setminus\{p,q\}$ the metric is preserved by $(A,A^{-1})$.
\end{proof}
Clearly $U$ is included in $U_A$,
hence $D$ defines a Lorentzian foliation on $\Sigma_2\times {\mathbb S}^1$. We denote it by $\mathcal F$. Let us remark that $\FF$ has the following interesting properties: it is Lorentzian, transversely Einstein and radial (i.e. its holonomy fixes a point). 

Let us tell now how to make more complicated examples from this one. First of all we can take $\hat \Pi : \widehat \Sigma \rightarrow \Sigma_2$ a finite cover and compose $D$ by $\hat \Pi$. But we can do better, let us take $\theta$ a closed $1$-form on $\widehat \Sigma $. Its pull-back $\widetilde\Pi^*\theta$ on the universal cover $\widetilde \Sigma$ is exact and therefore  there exists a function  $\Theta$ such that it is equal to $\dd \Theta$. We can replace the map $D$ by $D_\theta$ defined by 
 $$
\begin{array}{lccc}
& \widetilde \Sigma \times \R &\longrightarrow& \mathrm{Ein}_2\\
D_\theta:     &  (x,t) &\mapsto  &\left(\exp((t+\Theta(x))a)(P_1 \circ \widetilde \Pi (x)), \exp(-(t+\Theta(x))a)(P_2\circ \widetilde\Pi(x))\right).
\end{array}
$$
Those  maps define a whole family of radial, transversely Einstein, Lorentzian foliations $\FF_\theta$ on $\widehat \Sigma$. They depend on the cohomology class of $\theta$. Indeed, if  $\theta'=\theta+df$ with $f$ a function of $\widehat \Sigma$, the change of variable $t'=t+f(x)$ conjugates the foliations. It will be shown in a further work that this family contains a lot of examples not conjugated with each other.
\vskip .2cm
Let us describe now the dynamics of $\FF$. On figure 3, we draw the picture of $\FF$ restricted to $W\times {\mathbb S}^1$, where $W$ is a maximal open set of $\Sigma_2$ on which the projection $P$ is a diffeomorphism. The set $W$ is diffeomorphic to a disc with two holes and 
$\Sigma_2\times {\mathbb S}^1$ is actually given by two copies of $W\times {\mathbb S}^1$ glued along their boundaries $B_1\times{\mathbb S}^1$, $B_2\times{\mathbb S}^1$, $B_3\times{\mathbb S}^1$. We note that the $B_i's$ are circles. On figure 3, the vertical segments represent circles.
\begin{figure}[h]
\hspace*{\fill}
\epsfig{file=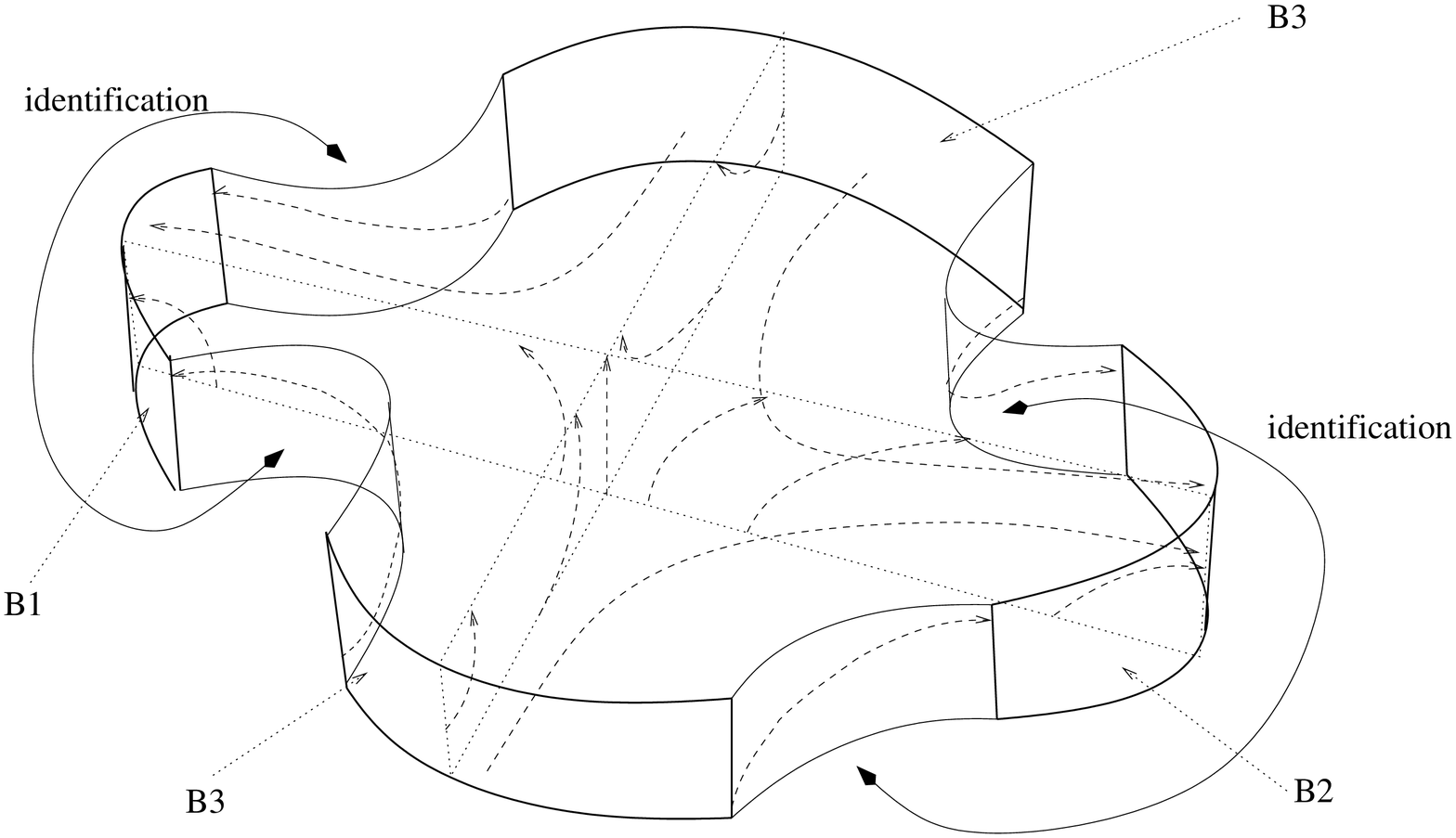,width=10cm}
\hspace*{\fill}\\
\hspace*{\fill}{\bf Figure 3.}\hspace*{\fill}
\end{figure}
We see that $\FF$ restricted to $W\times {\mathbb S}^1$ has only one closed leaf with nontrivial holonomy. 
On $\Sigma_2\times {\mathbb S}^1$ we have two interesting tori $T_1$ and $T_2$ given by the gluing of the dotted annuli of the picture, equal to $P^{-1}(\{A_+\}\times\R{\rm P}^1)\times{\mathbb S}^1$ and $P^{-1}(\R{\rm P}^1\times\{A_+\})\times{\mathbb S}^1$. Restricted to each of these tori $\FF$ is a Reeb foliation with two closed leaves. Finally, we can see that the leaves of $\FF$ not included in $T_1\cup T_2$ are compact. As equicontinuity is a delicate notion for foliations, admitting slightly different, and non-equivalent, definitions, we give the definition we consider.
\begin{defi}
Let $\FF$ be a $1$-dimensional foliation on a compact manifold $M$. Taking possibly a $2$-cover of $M$, the leaves of $\FF$ are the orbits of a flow $\varphi$ without fixed point. We say that $\FF$ is (uniformly) equicontinuous on an open $\FF$-saturated subset $U$ if $\{(\varphi_t)_{|U}, t\in\R\}$ is an (uniformly) equicontinuous set of diffeomorphisms of $U$.
\end{defi}

The foliation $\FF$ restricted to $(\Sigma_2\times {\mathbb S}^1)\setminus (T_1\cup T_2)$ is transversely parallelizable and so in particular equicontinuous. Moreover, $(\Sigma_2\times {\mathbb S}^1)\setminus (T_1\cup T_2)$ is an union of $\FF$-saturated open subsets, precompact in it, and on which $\FF$ is uniformly equicontinuous. Conversely on the $T_i$'s the foliation is nowhere uniformly equicontinuous. So the open dense subset $(\Sigma_2\times {\mathbb S}^1)\setminus(T_1\cup T_2)$ is the domain of equicontinuity of $\FF$. This kind of behaviour is not possible for an isometric flow (see for example section 3 of \cite{Zeghib-jdg}) therefore this foliation is not isometric. Actually the isometric Lorentzian foliations of $3$-manifolds are known, see \cite{Zeghib-jdg}, and are all transversely complete.  We do not give any proof that the dynamics of $\FF$ is as drawn on the picture but it can be  
checked by the reader.

To conclude we construct a codimension $1$ foliation transverse to $\FF$\label{exampletotogeo}. We consider first the codimension $1$ foliation $\mathcal G_0$ of $W\times{\mathbb S}^1$ whose leaves are the $W\times\{t\}$. This foliation is transverse to $\FF$ except along the $B_i\times{\mathbb S}^1$'s where it is everywhere tangent. Then we spin $\mathcal G_0$ along the three tori of the boundary (see \cite{cantwellconlon} p. 84 for example).  We do this spinning upwards along $B_3\times{\mathbb S}^1$ and downwards along $B_1\times{\mathbb S}^1$ and $B_2\times{\mathbb S}^1$ in order to remain transverse to $\FF$. We obtain a foliation $\mathcal G_1$ on $W\times {\mathbb S}^1$, tangent to the boundary, everywhere transverse to $\FF$, and whose holonomy along the toral leaves has a trivial Taylor development. Consequently it is possible to glue two copies of $(W\times {\mathbb S}^1, \mathcal G_1)$ along their boundary. The obtained foliation $\mathcal G$ is a smooth codimension $1$ foliation, transverse to $\FF$, with three toral leaves on $\Sigma_2\times {\mathbb S}^1$. According to proposition \ref{base} there exists a Lorentzian metric on $\Sigma_2\times{\mathbb S}^1$ for which $\mathcal G$ is timelike totally geodesic. This is the example announced at the end of section \ref{section classification 2-feuilletages}. As the orthogonal of $\GG$ is not transversely complete, its leaves are not geodesically complete. However it can be checked that the 3 compact leaves are geodesically complete because they are in the conformal class of a flat torus (see \cite{Ca-Ro}). A similar construction is possible for the foliations $\FF_\theta$. 
\subsection{A result without any assumption of completeness}
Up to a finite cover, a codimension $2$ Lorentzian foliation is given by the intersection of two codimension $1$ foliations, called its lightlike foliations. In the $3$ dimensional case, it is natural to wonder if those lightlike foliations may contain Reeb components. We prove the following.
\begin{prop}\label{Reeb}
Let $M$ be a closed $3$-manifold.
Let $\FF$ be a foliation of $M$ given by the intersection of two codimension $1$ foliations $\mathcal G_1$ and  $\mathcal G_2$.
If $\mathcal G_1$ contains a Reeb component then $\FF$ has no transverse volume form and therefore is not Lorentzian.
\end{prop}
\begin{proof}
This is a direct consequence of the results of I.\@ Tamura and A.\@ Sato \cite{tamurasato} about foliations transverse to a Reeb component. They proved that such a foliation must contain a half Reeb component. It is not hard to see that it implies that $\FF$ has an attractive leaf and therefore no transverse volume form. 
\end{proof}
Let us remark that this result can not be extended to taut foliations: the lightlike foliations of the transversely Einstein example given above are not taut: they have a separating toral leaf.
Anyway, thanks to the classical theorems of Novikov (see \cite{novikov}) and Palmeira (see \cite{palmeira}), we deduce the following theorem.
\begin{thm}
Let $M$ be a $3$-dimensional closed manifold. If $M$ possesses a Lorentzian foliation $\FF$ then its universal cover $\widetilde M$ is diffeomorphic to $\R^3$. Moreover the leaves of $\widetilde \FF$, the lift of $\FF$ to $\widetilde M$, are all lines.
\end{thm}
\begin{proof}
We know that the lightlike foliations of $\FF$ have no Reeb component then, by Novikov's theorem, we deduce that they have no vanishing cycles or equivalently that their lift to the universal cover as foliations by planes. Palmeira's theorem affirms then that the universal cover of $M$ is diffeomorphic to $\R^3$. To conclude we observe that $\widetilde \FF$ is tangent to a foliation by planes and therefore its leaves are lines.
\end{proof}

\section{\mathversion{bold}Codimension $2$ Lorentzian foliations and linear forms\mathversion{normal}}\label{formeslineaires}
It is well known that a transversely orientable codimension $1$ foliation $\FF$ is given by a $1$-form $\omega$ satisfying the Frobenius condition $\omega\wedge\dd\omega=0$. We know how to read the existence  of some transverse structure on $\FF$ on the $1$-form $\omega$. For example $\FF$ is transversely affine (resp.\@ projective)  if and only if there exists two $1$-forms $\omega$ and $\omega_0$ such that $\omega$ is an equation of the foliation with $\dd\omega=\omega\wedge\omega_0$ and $\dd\omega_0=0$ \cite{bobo} (resp.\@ there exists three $1$-forms $\omega$, $\omega_1$ and $\omega_2$ such that $\omega$ is an equation of the foliation with $\dd\omega=2\omega_1\wedge\omega$, $\dd\omega_1=\omega\wedge\omega_2$ and $\dd\omega_2=2\omega_2\wedge\omega_1$ \cite{blumenthal}). The problem with codimension $2$ foliations is that it is not true that up to a finite cover they are generated by a pair of $1$-forms. Their normal bundle must be trivial. But it is precisely the case, up to a possible $4$-cover, for codimension $2$ Lorentzian foliations. 

As we already saw, a codimension $2$ Lorentzian foliation preserves two normal line fields and a normal volume element.
Taking possibly a  $4$-cover this implies that the foliation is given by the intersection of two transversely orientable codimension $1$ foliations. Therefore there exist two $1$-forms $\omega_1$ and $\omega_2$ such that:
$$\left\{\begin{array}{c}
T\FF=\ker \omega_1 \cap \ker \omega_2\\
\dd \omega_i\wedge\omega_i=0,\ \mathrm{for}\  i=1,2
\end{array}\right.$$
But of course those conditions are not sufficient to imply that $\FF$ is Lorentzian. We did not use the fact that  $\FF$ admits a transverse volume form. This form can be written  $f\, \omega_1\wedge \omega_2$, with $f$ a positive function. But replacing $\omega_1$ by $f\, \omega_1$ we can suppose that the transverse volume form is $\omega_1\wedge\omega_2$. Thus, after a possible modification of the $\omega_i$, we have:
$$\mathcal L_X (\omega_1\wedge\omega_2)=0,$$
for any vector field $X$ tangent to $\FF$. We already have a characterization by the $1$-forms of being codimension $2$ Lorentzian.

\begin{prop}\label{3formes}
A codimension $2$ transversely orientable foliation $\FF$ is a transversely causal Lorentzian foliation if and only if there exists two  $1$-forms $\omega_1$ and $\omega_2$ such that for any vector field $X$ tangent to $\FF$ we have:
$$\left\{ \begin{array}{l}
T\FF=\ker\omega_1\cap\ker \omega_2,\\
\dd\omega_i\wedge \omega_i =0,\ \text{for }i=1,2\\
\mathcal{L}_X(\omega_1\wedge\omega_2)=0
\end{array}
\right.
$$
or, equivalently, there exists a third $1$-form $\omega_0$ such that: 
$$\left\{ \begin{array}{l}
T\FF=\ker\omega_1\cap\ker \omega_2,\\
\dd\omega_1=\omega_1\wedge\omega_0,\\
\dd\omega_2=-\omega_2\wedge\omega_0
\end{array}
\right.
$$
\end{prop}
\begin{proof}
The first assertion has been already proven. For the second one we just notice that $ \dd\omega_i\wedge \omega_i =0, \text{for } i=1,2$ implies that there exist two $1$-forms $\omega_3$ and $\omega_4$ such that $\dd\omega_i=\omega_i\wedge (\omega_{i+2}+f_i \omega_i)$, for any function $f_i$. Moreover we know  that $\mathcal L_X (\omega_1\wedge\omega_2)=0$ this implies that $i_X(\omega_{3}+f_1 \omega_1+ \omega_{4}+f_2 \omega_2)=0$, for any vectorfield $X$ tangent to $\FF$. We know moreover that there exists vectorfields $Y_1$ and $Y_2$ such that $\omega_i(Y_i)=0$ ($i=1$ or $2$) and $\omega_{3-i}(Y_i)$ is nowhere vanishing. Those vectorfields allow us to find the $f_i$'s such that 
$\omega_{3}+f_1 \omega_1= - (\omega_{4}+f_2 \omega_2)$. We denote this form by $\omega_0$.  
\end{proof}

\begin{rema}{\rm
We can give an alternative proof in the context of transverse coframes of the foliation. Let $\FF$ be a transversely causal Lorentzian foliation and let  $\GL^*(M,\FF)$ be the principal bundle of transverse coframes of $\FF$. The transverse structure gives a foliated reduction of $\GL^*(M,\FF)$ to the group ${\rm SO}_0(1,1)$. Let $\Fr(\FF)^*$ be this reduction. There exists a torsion free basic connection on $\Fr(\FF)^*$: that of Levi-Civita. It is a $\so(1,1)$-valued $1$-form ($\so(1,1)$ denotes the Lie algebra of $\SO(1,1)$), hence there exists a $1$-form on $\Fr(\FF)^*$: $\lambda$, such that the matrix of the connection is $\left(\begin{array}{cc}\lambda&0\\0&-\lambda\end{array}\right)$. This connection is torsion free i.e.\@ by definition $\dd\theta-\left(\begin{array}{cc}\lambda&0\\0&-\lambda\end{array}\right)\wedge\theta=0$ where the fundamental $1$-form $\theta$  of $\pi:\GL^*(M,\FF)\to M$ is given by $\theta_z(Z)=z(\overline{\dd_z\pi(Z)})$ (for $X\to \overline X$ the natural projection from $TM$ to $\nu(\FF)$). Since there exist two $1$-forms of $M$, say $\omega_1$ and $\omega_2$, such that $T\FF=\ker\omega_1\cap\ker \omega_2$, we can set $\bar\omega_i(\overline X)=\omega_i(X)$, for all vector $X$ in $T_xM$ and for $i=1,2$, and construct the (non basic) section $s|_x=(\bar\omega_1|_x,\bar\omega_2|_x)$ of $\Fr(\FF)^*$. As $s^*\theta|_x(X)=s|_x(\overline X)=(\omega_1|_x(X),\omega_2|_x(X))$ for all $x\in M$ and all $X\in TM$, the nullity of the torsion amounts to the equalities $\dd\omega_i\mp s^*\lambda\wedge\omega_i=0$ (for $i=1$ resp. $2$). Therefore the $1$-form $\omega_0$ of the last proposition is equal to the $1$-form $-s^*\lambda$.}
\end{rema}

This proposition enables us to give examples of codimension $2$ Lorentzian foliations. The case of algebraic Anosov flows of closed $3$-manifolds is perhaps not the most interesting. The $1$-forms are just the given by the dual  of the classical left invariant frames. We let the reader check the details.

The following is perhaps more interesting. On $\R^3$, we consider the forms $\omega_1=\cos x \dd x + \sin x \dd z$ and $\omega_2=\cos y \dd y -\sin y \dd z$. When both $x$ and $y$ are not equal to $\pi/2 \mod \pi$, this defines a codimension $2$ foliation. It is not hard to see that it gives a foliation of $\R^3/((2\pi\Z)^2\times\Z)\simeq\T^3$ minus the four circles given by $\cos x=\cos y=0$.
We have $\dd \omega_1 = \cos x \dd x\wedge \dd z= \omega_1 \wedge \dd z$ and $\dd \omega_2 = - \cos y \dd y\wedge \dd z=- \omega_2 \wedge \dd z$. Hence the foliation is Lorentzian. The lightlike foliations are just given by the product of a Reeb foliation of the $2$-torus by ${\mathbb S}^1$. Thus we can make the drawing of figure 4. The vertical faces of the cube are pieces of lightlike compact leaves, the dotted curves delimit pieces of lightlike leaves and the black curves are pieces of some leaves of the Lorentzian foliation given by the intersection of the lightlike leaves.
\begin{figure}[h]
\hspace*{\fill}
\epsfig{file=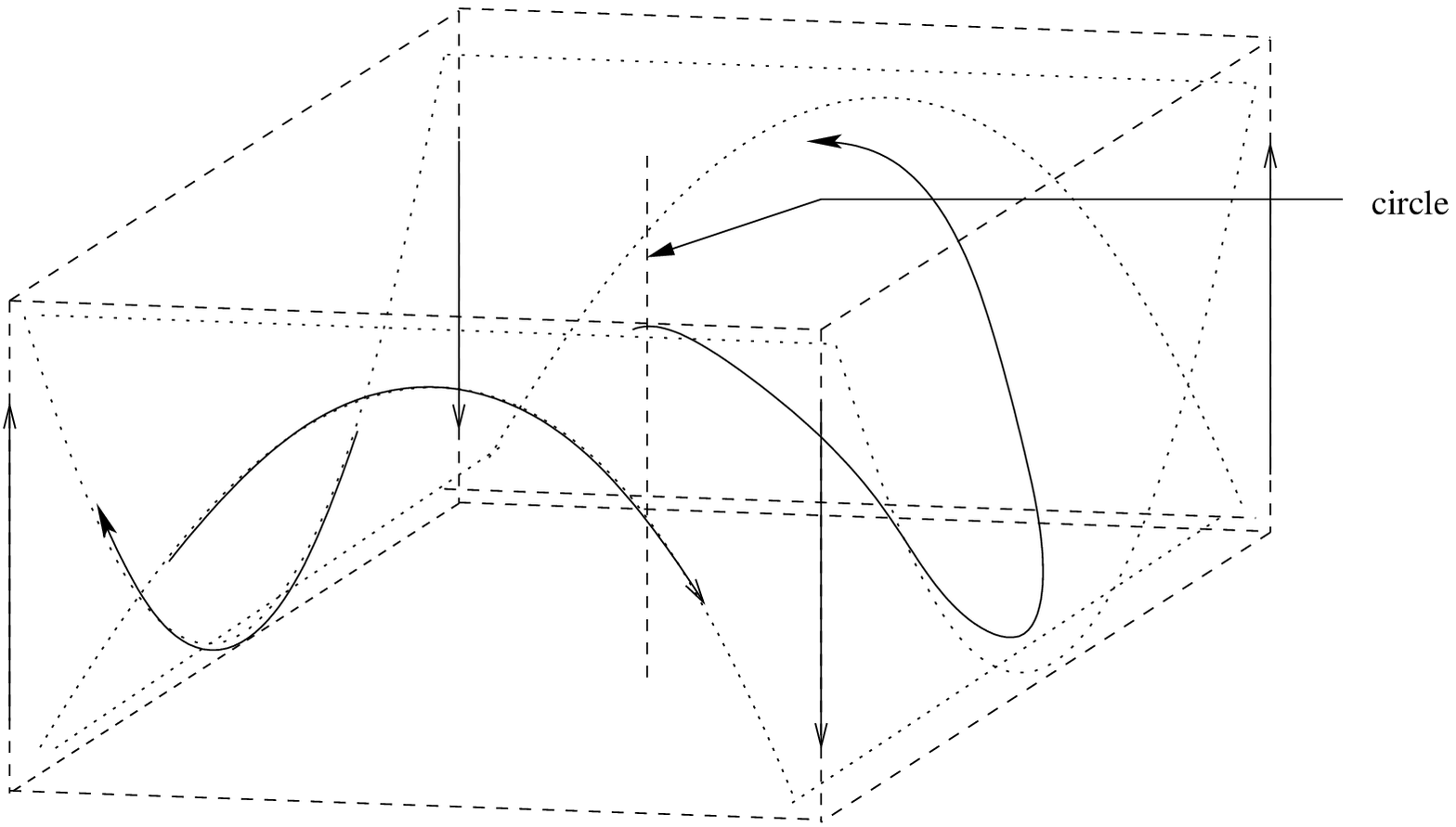,width=12cm}
\hspace*{\fill}\\
\hspace*{\fill}{\bf Figure 4.} The foliated parallelepiped $\left([0,\pi]^2\setminus\{(\frac{\pi}{2},\frac{\pi}{2})\}\right)\times{\mathbb S}^1$ in $\R^2\times{\mathbb S}^1$, \hspace*{\fill}\\
\hspace*{\fill}or in $(\R^2/(2\pi\Z)^2)\times{\mathbb S}^1$. The vertical segments represent the circles of the factor ${\mathbb S}^1$.\hspace*{\fill}
\end{figure}
This may remember something to the reader: this foliation, up to a four-cover, is diffeomorphic to that obtained in the beginning of section \ref{noncompletesection}, restricted to the open subset $P^{-1}((\R{\rm P}^1\setminus\{A_-\})^2)\times{\mathbb S}^1$ of $\Sigma\times{\mathbb S}^1$. Indeed, Corollary \ref{mink-sit} tells us that it is transversely Lorentzian flat (the form $\omega_0$ is closed).

\label{metriquenonplate} Another interesting property can be illustrated thanks to $1$-forms. We have seen that the Lie affine flow on $T^3_A$ is Lorentzian but that its space of transverse metrics is big. We are going to give an example of a transverse metric which is not the usual flat one. On $\R^3$ endowed with coordinates $(x,y,z)$, we consider the following $1$-forms:
$$\begin{array}{l}
\omega_1=\cos y \dd x + \lambda^{y/2\pi} \sin y \dd y,\\
\omega_2=\sin y \dd x - \lambda^{y/2\pi} \cos y \dd y,
\end{array}
$$
where $\lambda$ is a quadratic number such that $\lambda+1/\lambda \in \Z$.
Clearly, $\ker\omega_1\cap\ker \omega_2=\partial z$ and we have $\mathcal L_{\partial z}(\omega_1\wedge\omega_2)=0$. 
This gives a very particular transverse Lorentzian metric to this trivial foliation. Moreover those forms are invariant under any  translation fixing $y$, and by the transformation $\gamma_0 :\ (x,y,z)\mapsto (\lambda x, y + 2\pi, z/\lambda)$. There exist two translations $\gamma_1$ and $\gamma_2$ such that the group generated by $\gamma_0$,
$\gamma_1$ and $\gamma_2$  acts properly and cocompactly and preserves those $1$-forms. The quotient is just $T^3_A$ endowed with a foliation conjugated with the Lie affine foliation but the Lorentzian transverse structure induced by the forms has not constant curvature.
\vskip .2cm
Actually the form $\omega_0$ contains more geometric information about the transverse Lorent\-zian metric.  It follows classically from $\dd\omega_1=\pm\omega_i\wedge\omega_0$, $i=1,2$, that $0=\dd^2\omega_i=\pm(\dd\omega_i\wedge\omega_0-\omega_1\wedge\dd\omega_0)=\mp\omega_i\wedge\dd\omega_0$, so that there are two $1$-forms $\alpha$ and $\beta$ such that $\dd\omega_0=\alpha\wedge\omega_1$ and $\dd\omega_0=\beta\wedge\omega_2$. Finally there is a function $f$ such that $\dd\omega_0=f\omega_1\wedge\omega_2$. As $i_X\dd\omega_0=0$ and as $\dd\omega_0$ is closed, $\dd\omega_0$ is basic, thus $f$ is also basic. Besides $f$ depends only on the transverse metric, so it should be connected to the curvature of this metric. And indeed:
\begin{prop}\label{courb}
Let $\FF$ be a transversely orientable and causal Lorentzian foliation on a manifold $M$, let $K$ be the function on $M$ given by the transverse curvature  and let $\omega_0$, $\omega_1$, $\omega_2$ the $1$-forms given by proposition \ref{3formes}. Then 
$$\dd\omega_0= K \, \omega_1\wedge\omega_2.$$
\end{prop}
\begin{proof}
We simply recall the link, in dimension two, between the curvature form and the scalar curvature. The result follows. Let us denote by $\omega=\left(\begin{array}{cc}\lambda&0\\0&-\lambda\end{array}\right)$ the connection form of $\Fr(\FF)^\ast$; the induced curvature form is $\Omega=\dd\omega+\omega\wedge\omega=\dd\omega$ as, $\mathfrak{so}(1,1)$ being abelian, $\omega\wedge\omega$ vanishes. So $s^\ast\Omega={\rm diag}(s^\ast\dd\lambda,-s^\ast\dd\lambda)={\rm diag}(\dd(s^\ast\lambda),-\dd(s^\ast\lambda))={\rm diag}(\dd\omega_0,-\dd\omega_0)\in\bigwedge^2T^\ast M\otimes\mathfrak{so}(1,1)$ is basic as $\dd\omega_0$ is. So it may be viewed as an element of $\bigwedge^2\nu(\FF)^\ast\otimes\mathfrak{so}(1,1)$, which is a line bundle over $M$, identified with ${\rm End}(\bigwedge^2\nu(\FF)^\ast)$. Indeed, via the Lorentzian metric, $\bigwedge^2\nu(\FF)^\ast$ is identified with $\mathfrak{so}(1,1)$ by $\varphi:a\wedge b\mapsto \frac{1}{2}(a\otimes b^\sharp-b\otimes a^\sharp)$. So, this bundle has a canonical section, the constant ${\rm Id}_{\bigwedge^2\nu(\FF)^\ast}$, and $s^\ast\Omega=\sigma{\rm Id}$. This $\sigma$ is the curvature. Now, $\varphi(\bar{\omega}_1\wedge\bar{\omega}_2)=\left(\begin{array}{cc}1&0\\0&-1\end{array}\right)$, so $s^\ast\Omega=\sigma\left(\begin{array}{cc}\bar{\omega}_1\wedge\bar{\omega}_2&0\\0&-\bar{\omega}_1\wedge\bar{\omega}_2\end{array}\right)$. We are done.

One may also convince the reader by a calculation in transverse coordinates. On any small enough transversal $T$ to the foliation, endowed with the metric $g=\frac{1}{2}(\omega_1\otimes\omega_2+\omega_2\otimes\omega_1)$ defined by the pair of forms $\{\omega_1,\omega_2\}$, there are some coordinate systems $(x_1,x_2)$ such that the lines $\{x_i={\rm const.}\}$ are the integral leaves of $\omega_i$. Then $\omega_i=e^{f_i}\dd x_i$ with some functions $f_i$, and $g=e^f\dd x_1\dd x_2$ and $\omega_1\wedge\omega_2=e^f\dd x_1\wedge\dd x_2$, with $f=f_1+f_2$. We denote $\frac{\partial}{\partial x_i}$ by $X_i$. As ${\rm Span}(X_1)$ and ${\rm Span}(X_2)$ are the two isotropic distributions of $g$, they are parallel for its Levi-Civita connection $D$. So for all $i,j$: $D_{X_i}X_j\in{\rm Span}(X_j)$, hence, $D_{X_i}X_j=0$ for $i\neq j$ and $D_{X_i}X_i=\frac{\partial f}{\partial x_i}X_i$. In such coordinates, it is noticeable that the curvature $\sigma$ satisfies:
\begin{align*}
\sigma&=g(R(X_1,X_2)X_1,X_2)/g(X_1,X_2)^2\\
&=e^{-2f}g(D_{X_1}D_{X_2}X_1-D_{X_2}D_{X_1}X_1,X_2)\\
&=e^{-2f}(-(\partial g(D_{X_1}X_1,X_2)/\partial x_2)+g(D_{X_1}X_1,D_{X_2}X_2))\qquad\text{as $D_{X_2}X_1\equiv 0$,}\\
&=e^{-2f}(-(\partial ((\partial f/\partial x_1)e^f)/\partial x_2)+g(({\partial f}/{\partial x_1})X_1,({\partial f}/{\partial x_2})X_2))\\
&=e^{-2f}(-(\partial^2f/\partial x_1\partial x_2)e^f-(\partial f/\partial x_1)(\partial f/\partial x_2)e^f+(\partial f/\partial x_1)(\partial f/\partial x_2)e^f)\\
&=-(\partial^2f/\partial x_1\partial x_2)e^{-f}
\end{align*}
On their side, $\dd\omega_i=\frac{\partial f_i}{\partial x_{3-i}}e^{f_i}\dd x_{3-i}\wedge\dd x_i=(-1)^{3-i}\omega_i\wedge\omega_0$, with $\omega_0=\frac{\partial f_2}{\partial x_1}\dd x_1-\frac{\partial f_1}{\partial x_2}\dd x_2$. So $\dd\omega_0=-(\frac{\partial^2 f_1}{\partial x_1\partial x_2}+\frac{\partial^2 f_2}{\partial x_2\partial x_1})\dd x_1\wedge\dd x_2=-\frac{\partial^2f}{\partial x_1\partial x_2}e^{-f}\omega_1\wedge\omega_2$, so the result.
\end{proof}
Thanks to this proposition we can find the Lorentzian analogue to Blumenthal's result (see \cite{blumenthal}) about codimension $2$ Riemannian foliations with transverse constant curvature  (this could also be done by Blumenthal's method).
We have a characterization of codimension $2$ de Sitter and Minkowski foliations.
\begin{cor}\label{mink-sit}
Up to a possible $4$-cover, a codimension $2$ foliation is transversely Minkowski (resp.\@ de Sitter) if and only if there exists three $1$-forms $\omega_0$, $\omega_1$ and $\omega_2$ such that
$$\left\{ \begin{array}{l}
T\FF=\ker\omega_1\cap\ker \omega_2,\\
\dd\omega_1=\omega_1\wedge\omega_0,\\
\dd\omega_2=-\omega_2\wedge\omega_0\\
\dd \omega_0 = 0\quad (\mathrm{ resp.}\ \omega_1 \wedge \omega_2)
\end{array}
\right.
$$
\end{cor}
\begin{rema}
 {\rm If we consider an orientable, time-orientable,  Lorentzian surface $(S,g)$, it  has a codimension $2$ Lorentzian foliation by  points and proposition \ref{courb} is still true. It says that if $S$ is compact the integral 
 of $K\, \mathrm{vol}_g$ vanishes. As in this case $S$ is diffeomorphic to a $2$ torus, this is  
  the Gauss-Bonnet-Avez theorem. This also shows that the canonical volume form of a locally de Sitter surface $S$ is exact and therefore, as it is well known, $S$ is certainly not compact. }
 \end{rema}

\noindent 
C.\@ Boubel,\hfill{\tt Charles.Boubel{\rm @}ens-lyon.fr}
\vskip .1cm\noindent 
C.\@ Tarquini,\hfill{\tt Cedric.Tarquini{\rm @}ens-lyon.fr}
\vskip .2cm\noindent 
{\'Ecole Normale Sup\'erieure de Lyon,\\
 Unit\'e de Math\'ematiques Pures et Appliqu\'ees,\\
  46 all\'ee d'Italie,\\
 F-69364 LYON CEDEX 07\\
 
\noindent 
 P.\@ Mounoud,\hfill {\tt pierre.mounoud{\rm @}math.u-bordeaux1.fr}
 \vskip .2cm\noindent 
 Universit\'e Bordeaux 1,\\ Laboratoire Bordelais d'Analyse et G\'eom\'etrie,\\
  351 cours de la Libération,\\
   F-33405 TALENCE CEDEX\\

\end{document}